\documentclass[review]{elsarticle}
\usepackage{lineno,hyperref}
\usepackage{graphicx}
\usepackage{amsthm}
\usepackage{amssymb}
\usepackage{amsmath}
\usepackage{bbm}
\usepackage{caption}
\usepackage{subcaption}
\usepackage{amsfonts}
\usepackage{xcolor}
\usepackage{lineno,hyperref}
\usepackage{epstopdf}
\newtheorem{proposition}{Proposition}

\newtheorem{remark}{Remark}

\modulolinenumbers[5]
\begin{document}

\markboth{B. Aymard}{Oscillations, chaos and strange attractors in presence of self-cross diffusion}

\title{Oscillating Turing patterns, chaos and strange attractors 
in a reaction-diffusion system 
augmented with self- and cross-diffusion terms}

\author{Benjamin Aymard}

\address{MathNeuro Team,\\ 
Inria branch at the University of Montpellier,\\
Montpellier, France}

\begin{abstract}
In this article we introduce an original model in order to study the emergence of chaos in a reaction diffusion system in the presence of self- and cross-diffusion terms.
A Fourier Spectral Method is derived to approximate equilibria and orbits of the latter.
Special attention is paid to accuracy, a necessary condition when one wants to catch periodic orbits and to perform their linear stability analysis via Floquet multipliers.
Bifurcations with respect to a single control parameter are studied in four different regimes of diffusion:
linear diffusion, self-diffusion for each of the two species, and cross-diffusion.
Key observations are made: 
development of original Turing patterns, 
Hopf bifurcations leading to oscillating patterns 
and period doubling cascades leading to chaos.
Eventually, original strange attractors are reported in phase space. 
\end{abstract}

\maketitle


\section{Introduction}

%Chaos theory
When does a deterministic dynamical system become unpredictible?
In his 1814 memoir on probabilities \cite{Laplace:1814}, Pierre-Simon de Laplace hypothetised that one could predict the future of a dynamical system given the exact knowledge of positions and impulsion of its components at an initial time.
Later, James Clerk Maxwell observed in 1876 \cite{Maxwell:1876} that for a class of systems, an error on the initial condition could dramatically modify its outcome. 
Eventually, Henri Poincaré's work in celestial mechanics from 1892 \cite{Poincare:1892} led him to discover the unpredictability of trajectories in the three body problem.
In doing so, he provided fundamental tools, of what would later become the chaos theory. 
Beyond its philosophical aspect, and its origin in celestial mechanics, chaos theory finds practical applications in everyday life, in problems ranging from meteorology \cite{Lorenz:1963} to electrical engineering, passing by chemical reactions \cite{Rossler:1976}, plasma physics \cite{Chirikov:1959}, laser physics \cite{Popov:1992} and mechanics \cite{Yorke:1991}. 
Understanding chaotic systems is a challenging task, with numerous potential applications.

%bifurcation and chaos
Using the tools of bifurcation theory \cite{Manneville:2004, Kuznetsov:2004}, 
several mechanisms explaining the developement of chaos have been identified.
The key idea of this aproach is to evaluate the sequence of bifurcations
(loss of stability of equilibria or periodic orbits)
while a control parameter is varied. 
%period doubling cascade
Among the identified scenarios is the period doubling cascade. 
Discovered in 1976 in numerical experiments on the logistic map 
(see \cite{May:1976}),
then observed in numerous experiments, such as in electronic circuits in 1981 by Linsay \cite{Linsay:1981}, or in mercury convection rolls in 1982 by Libchaber \cite{Libchaber:1982},  
it consists of a sequence of period doubling of periodic orbits. 
After each period doubling, the distance between consecutive bifurcation diminishes,
and eventually, the period tends to infinity for a finite value of the control param€ter: 
the orbit becomes aperiodic.
Careful measurements of distances between consecutive bifurcation values led to the discovery of the Feigenbaum constant \cite{Feigenbaum:1978, Coullet:1978}, revealing a concept of universality of the period doubling cascade.
It was then observed in many other models, in particular in the 
R\"{o}ssler model for chemical reaction \cite{Rossler:1976}, 
the Mackey-Glass system in blood cells population dynamics \cite{MackeyGlass:1977},
neuron models of Fitzhugh-Nagumo type \cite{Ermentrout:1984}, 
the Kuramoto-Sivashinsky system modeling flamme front propagation \cite{Papageorgiou:1991}, 
the Chua circuit \cite{Chua:1992}, the forced Dufing oscillator \cite{Duffin:1918}, among others.
Interested reader may find a historical note about period doubling in \cite{Collet:2019}.

%Roads and crossroads
However, even if many road to chaos have been identified, 
the question of understanding how does the presence of self- and/or cross-diffusion affect the development of chaos in reaction diffusion systems, remains unclear.
This work, which aims precisely at answering this question, 
stands at the crossroads of three domains that appeared during the 20th century.
First, the domain of morphogenesis, initiated by the seminal paper of Alan Turing \cite{Turing:1952}, studying the appearance of periodically spatial patterns, as a result of an instability created by the interaction between linear diffusion, and nonlinear reaction.
Second, the domain of spatio-temporal chaos, notably following the works of Kuramoto on diffusion induced chaos, where, for instance, the interaction between a Hopf bifurcation and a Turing instability 
can create oscillating patterns, possibly becoming chaotic \cite{Hopf:1942,Kuramoto:1978}.
Third, the domain of self- and cross-diffusion, motivated by questions in population dynamics and in chemistry, 
notably the works of Shigesada Kawasaki Terramoto \cite{SKT:1979} and Vanag and Epstein \cite{Epstein:2003}, and more recently in \cite{Moussa:2019, Aymard:2023}.

%Challenges
There is therefore a need for a simple model, with bifurcations driven by a single parameter, 
gathering those three elements, for which we could study in a continuous way the development of chaos.
Due to the presence of nonlinearities on both reaction and diffusion terms, this kind of system is, 
in general, not analytically solvable, except in special cases. 
Therefore, reliable numerical methods \cite{Tuckerman:2004} have to be designed in order to study them numerically. 
To this end, high accuracy is a necessary condition when one wants to catch periodic orbits, 
and to perform their linear stability analysis via Floquet multipliers \cite{Floquet:1883}.

%Outline and results
The outline of the paper is the following.
First, an original model is introduced, containing self- and cross-diffusion, 
with bifurcations controled by a single parameter. 
Together with the model, a Fourier spectral method is derived, 
and theoretical convergence guarantees are provided.
Second, bifurcations, of both equilibria and periodic orbits, are studied in four regimes: linear diffusion, activator self-diffusion, inhibitor self-diffusion, and cross-diffusion.  
Finally, the road to chaos is identified in each case, and the strange attractors that appear in the chaotic regime, are described.

\section{Model and method}

\subsection{Model}

In \cite{BVAM:1999, BVAM:2012} the authors introduce a model containing a richness of behavior, while being relatively simple, as the bifurcations are driven by one control parameter $C$, making it, as they mention, a good laboratory for studying oscillation and chaos in two species reaction diffusion models.
We will refer to this model as the classical BVAM model.
Along those lines, our aim is to introduce a generalization of their model, by introducing non linear diffusion terms.
Let us consider the dynamics of two species, $u_1$ and $u_2$, on a domain $\Omega = [-L_x,L_x]$, described by the reaction diffusion system:
\begin{align}
\frac{\partial u_1}{\partial t} &= 
\Delta (d_1 u_1 + d_{11}u_1^3 + d_{12}u_2^2u_1)
+ \eta(u_1 + a u_2 - C u_1 u_2 - u_1u_2^2),\nonumber\\
\frac{\partial u_2}{\partial t} &= 
\Delta (d_2 u_2 + d_{22}u_2^3 + d_{12}u_1^2u_2)
+ \eta( b u_2 + H u_1 + C u_1 u_2 + u_1u_2^2),
\label{model}
\end{align}
with $\eta$, $a$, $b$, $C$ and $H$ reaction parameters, and $d_i, d_{ii}, d_{ij}$ for $i,j=1,2,i\not=j$ diffusion parameters. Let us remark that, when $d_{11} = d_{22} = d_{12} = 0$, the model (\ref{model}) reduces to the original BVAM model.
The problem is mathematically closed by adding initial conditions and periodic boundary conditions on the domain $\Omega$. 
In their original work \cite{BVAM:2012}, the authors considered Dirichlet boundary conditions. 
Our choice of periodic boundary conditions here is motivated by two reasons:
this reduces the importance of the central point of the domain, later used to define a phase space, 
and it allows to easily use Fourier analysis \cite{GW:1999}.
From \cite{Aymard:2022} we know that models of the form (\ref{model}) follow an energy law of the form:
\begin{equation}
\frac{d}{dt}E(t) 
= 
-\sum_{i=1}^2 \|\nabla \mu_i \|_{L^2(\Omega)}^2 
+ \sum_{i=1}^2(R_i,\mu_i),
\label{dE}
\end{equation}
with $E$ defined by:
\begin{equation}
E(t) = \int_{\Omega} \left( d_1 \frac{u_1^2}{2} + d_2 \frac{u_2^2}{2} + d_{11} \frac{u_1^4}{4} + d_{22} \frac{u_2^4}{4} + d_{12} \frac{u_1^2u_2^2}{2} \right)dx,
\label{E}
\end{equation}
with reaction terms $R_1,R_2$ defined by:
\begin{align}
R_1 &= \eta(u_1 + a u_2 - C u_1 u_2 - u_1u_2^2), \nonumber\\
R_2 &= \eta( b u_2 + H u_1 + C u_1 u_2 + u_1u_2^2), \label{R}
\end{align}
and chemical potentials $\mu_1, \mu_2$ defined by:
\begin{align}
\mu_1 &= d_1 u_1 + d_{11}u_1^3 + d_{12}u_2^2u_1, \nonumber\\
\mu_2 &= d_2 u_2 + d_{22}u_2^3 + d_{12}u_1^2u_2. \label{mu}
\end{align}
Due to the complexity of the model, analytical solutions may only exist in special cases;
therefore, an efficient numerical method is needed in order to approximate orbits of this dynamical system.

\subsection{Fourier Spectral Method}

\subsubsection{Spatial approximation}

Using periodicity assumption, one may consider the Fourier series (in space):
\begin{align}
u_i(x,t) &= \sum_{k=-\infty}^{\infty} (\hat{u_i})_k(t) e^{i k x / |\Omega|},\quad i=1,2\\
\mu_i(x,t) &= \sum_{k=-\infty}^{\infty} (\hat{\mu_i})_k(t) e^{i k x / |\Omega|},\quad i=1,2\\
R_i(x,t) &= \sum_{k=-\infty}^{\infty} (\hat{R_i})_k(t) e^{i k x / |\Omega|},\quad i=1,2,
\label{FS}
\end{align}
with Fourier coefficients defined by:
\begin{align}
(\hat{u_i})_k(t) &= \frac{1}{|\Omega|}\int_{\Omega}u(x,t) e^{-i k x / |\Omega|}  \mathrm{d}x,\quad i=1,2\\
(\hat{\mu_i})_k(t) &= \frac{1}{|\Omega|}\int_{\Omega}\mu(x,t) e^{-i k x / |\Omega|} \mathrm{d}x,\quad i=1,2\\
(\hat{R_i})_k(t) &= \frac{1}{|\Omega|}\int_{\Omega}R(x,t) e^{-i k x / |\Omega|} \mathrm{d}x, \quad i=1,2.
\label{Fc}
\end{align}
Injecting those expressions in the model (\ref{model}), and using (\ref{R}), (\ref{mu}), one may write:
\begin{equation}
\frac{\partial \hat{u}_i}{\partial t} = -k^2 \hat{\mu}_i + \hat{R}_i, \quad i=1,2.
\label{FT_EDP}
\end{equation}
In Fourier space, the set of partial differential equations has become a set of ordinary differential equations, 
easier to solve.
However, except in simple cases, the Fourier series have to be evaluated numerically.
To this end, we use the Fast Fourier Transform algorithm \cite{GW:1999}.
Given a discretization of the domain $\Omega = [-L_x,L_x]$ by a step: 
\[
\Delta x = \frac{2 L_x}{N},
\]
leading to collocation points:
\[
x_j = j \Delta x \mbox{ with } \quad j= -\frac{N}{2},...,\frac{N-1}{2}.
\]
we compute an approximation of (\ref{Fc}), called the discrete Fourier transform:
\[
\tilde{\mu}_k = \frac{1}{N}\sum_{j=0}^{N-1} \mu_j e^{-i k x_j},
\]
Frequencies $f_m$ are evaluated as:
\[
f_m = m \Delta f 
\]
with:
\[
\Delta f = \frac{1}{N \Delta x} = \frac{1}{2 L_x}.
\]
Eventually, wave numbers $k$ and frequencies $f$ are related as:
\[
k = 2 \pi f.
\]

\begin{remark}
Several conventions coexist for the FFT algorithm.
The results are independent from the used convention, however,  
it is crucial to respect the FFT convention of the used library when one defines $k$, in order to be consistent during all the FFT and inverse FFT operations.
\end{remark}

\begin{remark}
The reader may note that nonlinear terms are first evaluated in real space, then transformed in Fourier space.
\end{remark}

The spatial approximation verifies an important convergence property, that we will now demonstrate.
Let us first recall the definition of a type of Sobolev space \cite{GW:1999}:
\[
H^{m} = \{f \in L^2, (1+|k|^2)^{m/2}\hat{f} \in L^2 \}
\]

\begin{proposition}
Let us assume that model (\ref{model}) admits equilibrium solutions:
\begin{equation}
-\Delta \mu_i(u_1,u_2) = R_i(u_1,u_2).
\label{eq}
\end{equation}
Let us further assume that those solutions are bounded:
\[
\| u_i \|_{L^2} \leq C.
\]
and smooth enough such that $u_i \in H^m$.
Then, for $N$ large enough, the following inequality holds:
\begin{equation}
\| \mu_i - P_N(\mu_i)\|^2_{L^2}  \leq \frac{\| R_i \|^2_{H^m} }{(1+N^{2})^m N^4}.
\label{convergence}
\end{equation}
\end{proposition}

\begin{remark}
The assumptions result from observations of numerous simulations with various parameters.
However, the proof of existence, regularity and boundedness of solutions remain open.
\end{remark}

\begin{proof}
By application of Fourier synthesis theorem \cite{GW:1999}, one gets:
\begin{align*}
\mu_i(x) &= \sum_{k=-\infty}^{\infty} \hat{\mu_i}_ke^{ikx}\\
&=
\underbrace{\sum_{k=-N/2}^{N/2} (\tilde{\mu_i})_ke^{ikx}}_{P_N(\mu)}
+ \underbrace{\sum_{k=-N/2}^{N/2} ((\hat{\mu_i})_k - (\tilde{\mu_i})_k)e^{ikx}}_{A_N(\mu)}
+ \underbrace{\sum_{|k| > N/2} (\hat{\mu_i})_ke^{ikx}}_{T_N(\mu)},
\end{align*}
with $P_N$ the interpolation polynom, $A_N$ the aliasing error, and $T_N$ the truncation error.
We then get, using orthogonality of the Fourier basis:
\[
\| \mu_i - P_N(\mu_i) \|^2_{L^2} = \| A_N(\mu_i) \|^2_{L^2} + \| T_N(\mu_i) \|^2_{L^2}.
\]
For the aliasing error $A_N$, let us recall the classical argument, relying on periodicity of the complex exponential, and the Fourier series evaluated at collocation points:
\[
\mu_i(x_j) 
= \sum_{k=-\infty}^{\infty} (\hat{\mu_i})_k e^{ikx} 
= \sum_{k=0}^{N-1} \left( \sum_{q=-\infty}^{\infty} (\hat{\mu_i})_{(k+qN)} \right) e^{ikx},
\]
then, by identification, one gets:
\[
(\hat{\mu_i})_k - (\tilde{\mu_i})_k 
= 
\sum_{q=-\infty, q\not=0}^{\infty} (\hat{\mu_i})_{(k+qN)}.
\]
For $N$ large enough, as $\mu_1,\mu_2$ are regular, the aliasing error is null.
For the truncation error, one may expand, using equation (\ref{FT_EDP}):
\[
\| T_N(\mu_i) \|^2_{L^2} 
= \sum_{|k| > N/2} |(\hat{\mu_i})_k|^2 
= \sum_{|k| > N/2} \frac{|(\hat{R_i})_k|^2}{k^4}.
\]
Then, multiplying numerator and denominator by the same constant, then by using the hypothesis on regularity, and finally the Parseval identity, one gets:
\begin{align*}
\| T_N(\mu_i) \|^2_{L^2} 
&= \sum_{|k| > N/2} \frac{(1+k^2)^m}{(1+k^2)^m}\frac{|\hat{R_i}_k|^2}{k^4}
\leq \frac{1}{(1+N^2)^m N^4} \sum_{|k| > N/2} (1+k^2)^m|\hat{R_i}_k|^2\\
&\leq \frac{\| R_i \|^2_{H^m} }{(1+N^{2})^mN^4}.
\end{align*}
As $u_1,u_2$ are bounded, by hypothesis, and as $R_1,R_2$ are polynomials in $u_1,u_2$, 
then for all $m$, $R_i \in H^m$. 
\end{proof}

\subsubsection{Temporal approximation}

To approximate the temporal part, the classical Runge-Kutta method of fourth order (RK4) is used (\cite{Demailly:2004})
Defining a time step $\Delta t$, and denoting by $\hat{u}_i^n$ the approximation of $\hat{u}_i$ at time $t^n = n \Delta t$, the RK4 scheme reads:
\begin{equation}
\hat{u}_i^{n+1} = \hat{u}_i^n + \frac{\Delta t}{6}\left(\hat{F}(U_1) + 2 \hat{F}(U_2) + 2\hat{F}(U_3) + \hat{F}(U_4) \right),
\label{RK4}
\end{equation}
with:
\begin{align*}
U_1 &= (\hat{u}_1,\hat{u}_2),\\ 
U_2 &= U_1 + \frac{\Delta t}{2}\hat{F}(U_1),\\ 
U_3 &= U_1 + \frac{\Delta t}{2}\hat{F}(U_2), \\
U_4 &= U_1 + \Delta t \hat{F}(U_3).
\end{align*}

\section{Patterns and oscillations}

We study the bifurcations of the system with respect to parameter $C$,
in the presence of four different kind of diffusion: a classical linear diffusion (seen as a witness test, the original BVAM model),
a self-diffusion on each of the two species, and a cross-diffusion between the two species. Parameters are gathered in \ref{parameters}.
Let us recall that in their original article, \cite{BVAM:2012} considered $\Delta x = 0.2$ and $\Delta t = 0.01$.

%\centering
\captionof{table}{Parameters}
\label{parameters}
\begin{tabular}{|c|c|}
\hline
Reaction & $L_x = 5 $, $H=3 $,$\eta = 1$, $a = -1$, $b = -3/2$ \\
\hline
Linear diffusion & $d_1 = 0.08$, $d_2 = 1 $, $d_{11} = 0, d_{12} = 0, d_{22} = 0$\\
Self-diffusion on $u_1$ & $d_1 = 0.08$, $d_2 = 1 $, $d_{11} = 0.07,d_{12} = 0 d_{22} = 0$\\
Self-diffusion on $u_2$ & $d_1 = 0.08$, $d_2 = 1 $, $d_{11} = 0, d_{12} = 0, d_{22} = 0.05$\\
Cross-diffusion & $d_1 = 0.08$, $d_2 = 1 $, $d_{11} = 0, d_{12} = 0.02, d_{22} = 0$\\
\hline
Spectral method (periodic orbits) & $N=300$ , $\Delta t = 8e^{-5}$\\
Spectral method (strange attractors) & $N=500 $, $\Delta t = 2e^{-5}$\\
Newton-Krylov method & $N=1000 $, $\epsilon = 1e^{-10}$\\
\hline
\end{tabular}
%, $n_1=n_2=2$, $A_1=0.8$, $A_2 = 0.4$

\subsection{Bifurcation of equilibria}

\subsubsection{Equilibria}

Equilibria of dynamical system (\ref{model}), denoted by $(\bar{u}_1,\bar{u}_2)$, 
are defined as stationnary solutions, verifying:
\begin{equation}
\frac{\partial}{\partial t}
\begin{pmatrix}
\bar{u_1}\\
\bar{u_2}
\end{pmatrix}
=
\begin{pmatrix}
\Delta \mu_1(\bar{u}_1,\bar{u}_2) + R_1(\bar{u}_1,\bar{u}_2)\\
\Delta \mu_2(\bar{u}_1,\bar{u}_2)  + R_2(\bar{u}_1,\bar{u}_2)
\end{pmatrix}
=
\begin{pmatrix}
0\\
0
\end{pmatrix}
\label{equilibrium}
\end{equation}
Homogeneous steady states are spatially constant equilibria. 
With the choice of parameters given in Table \ref{parameters}, 
BVAM have showed that only one solution exist, given by $(\bar{u}_1 = \bar{u}_2 = 0)$.
In the presence of sufficiently high diffusion, a Turing instability appears, creating non homogeneous states, 
that can not be analytically computed in general.
In order to evaluate equilibria (\ref{equilibrium}), we use the Newton-Krylov method \cite{NK:2004} to find the roots of the equivalent system:
\begin{equation}
\begin{cases}
\mathcal{F}^{-1}(\hat{F_1}(u_1,u_2)) = 0, \\
\mathcal{F}^{-1}(\hat{F_2}(u_1,u_2)) = 0. 
\end{cases}
\label{steady}
\end{equation}
until a given tolerance $\epsilon$ is reached, starting from an initial guess of the form:
\[
u_i(x) = A_i\cos(n_i x), i=1,2.
\] 
with $n_i$ a parameter, tuning the wave number of the initial guess.

\subsubsection{Stability of equilibria}

Linear stability analysis consists of studying the stability of equilibria, 
by studying the flow of linear perturbation around steady states. 
Given a perturbation $(\delta u_1,\delta u_2)$ around an equilibrium $(\bar{u_1},\bar{u_2})$, 
the linearized flow around a steady state reads:
\[
\frac{\partial}{\partial t}
\begin{pmatrix}
\delta u_1\\
\delta u_2
\end{pmatrix}
\approx
L(\bar{u}_1,\bar{u}_2)
\begin{pmatrix}
\delta u_1\\
\delta u_2
\end{pmatrix}
\]
with:
\begin{align}
\begin{split}
&L(\bar{u}_1,\bar{u}_2)
= \\
&\begin{pmatrix}
\frac{\partial R_1}{\partial u_1}(\bar{u}_1,\bar{u}_2)\delta u_1 + \Delta\left(\frac{\partial \mu_1}{\partial u_1}(\bar{u}_1,\bar{u}_2)\delta u_1\right) &
\frac{\partial R_1}{\partial u_2}(\bar{u}_1,\bar{u}_2)\delta u_2 + \Delta\left(\frac{\partial \mu_1}{\partial u_2}(\bar{u}_1,\bar{u}_2)\delta u_2\right)\\
\frac{\partial R_2}{\partial u_1}(\bar{u}_1,\bar{u}_2)\delta u_1 + \Delta\left(\frac{\partial \mu_2}{\partial u_1}(\bar{u}_1,\bar{u}_2)\delta u_1\right)&
\frac{\partial R_2}{\partial u_2}(\bar{u}_1,\bar{u}_2)\delta u_2 + \Delta\left(\frac{\partial \mu_2}{\partial u_2}(\bar{u}_1,\bar{u}_2)\delta u_2\right)
\end{pmatrix}
\end{split}
\label{L}
\end{align}
Using FFT, one may approximate operator $L$ in finite dimension. 
Indeed, as the Fourier Transform is linear, one can evaluate it on a basis, to get a matrix form:
\[
F_N = \hat{I}_n.
\]
Using this form allows to express Laplace operator as:
\[
\Delta_N = F_N^{-1} \times K \times F_N.
\]
with $K$ a diagonal matrix, such that $K_{mm} = k_m^2$ and $K_{mn} = 0$ for $m \not = n$. 
Replacing $\Delta$ operator in (\ref{L}) by $\Delta_N$ defines a matrix, 
from which stability is evaluated by computing its spectrum.

\subsubsection{Results}

\begin{figure}[!htbp]\centering
\begin{subfigure}{0.45\textwidth}
\includegraphics[width=\linewidth]{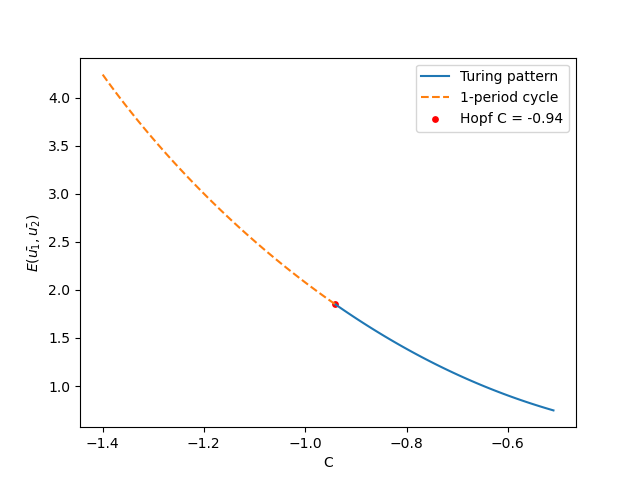}
\caption{Linear diffusion}
\end{subfigure}
\begin{subfigure}{0.45\textwidth}
\includegraphics[width=\linewidth]{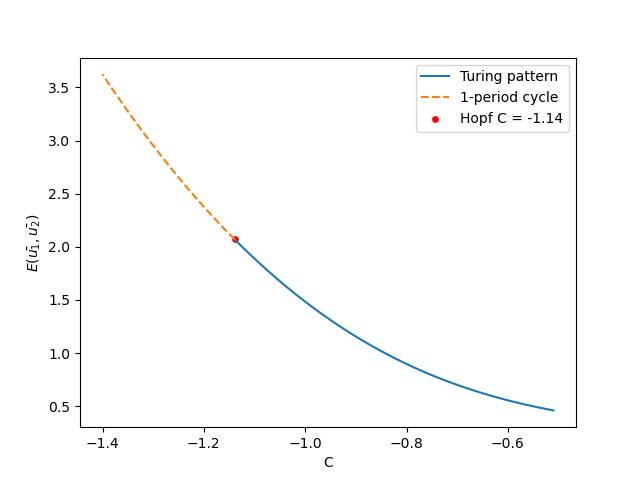}
\caption{Self-diffusion on $u_1$}
\end{subfigure}
\begin{subfigure}{0.45\textwidth}
\includegraphics[width=\linewidth]{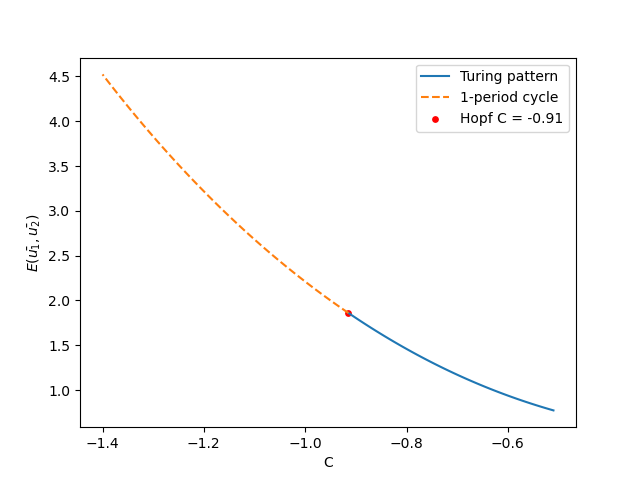}
\caption{Self-diffusion on $u_2$}
\end{subfigure}
\begin{subfigure}{0.45\textwidth}
\includegraphics[width=\linewidth]{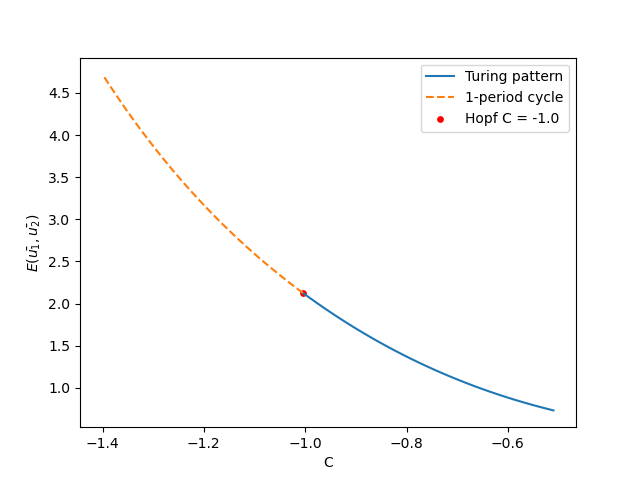}
\caption{Cross-diffusion}
\end{subfigure}
\caption{Bifurcation of equilibria. 
Steady states (\ref{equilibrium}) of model (\ref{model}), are computed using (\ref{steady}), and the corresponding energy (\ref{E}) is displayed, with respect to bifurcation parameter $C$.
A Hopf bifurcation appears in each case, when a stable equilibrium (continuous line), becomes unstable (dotted line), 
as a couple of complex eigenvalues with non negative real part appear in the spectrum (\ref{L}).
Past the Hopf bifurcation, the Turing patterns starts to oscillate in time, creating a periodic orbit.}
\label{Hopf}
\end{figure}

Using a parameter continuation strategy \cite{Kuznetsov:2004} on parameter $C$, 
equilibria and their stability have been computed in four regimes: linear diffusion, self-diffusion on $u_1$, then on $u_2$, and eventually cross-diffusion.
Results are displayed on Figure \ref{Hopf}. 
For $C$ varying from $-0.5$ to $-1.5$ in $100$ steps, steady states are computed, then the corresponding energy (\ref{E}) is evaluated. 
Evaluating the $L$ matrix (\ref{L}) on the steady states, stability of equilibria is also evaluated.
Let us note here that a threshold of $1e-3$ is considered on the real part of eigenvalue, to avoid numerical false positives.
For each of the four cases, the energy curve grows non linearly as parameter $C$ decreases.
Equilibrium is stable at first (continuous line), as the eigenvalues of $L$ are all with a strictly negative real part. 
Then, stability is lost (dotted line) through a Hopf bifurcation \cite{Hopf:1942} (indicated as a red dot), when a pair of complex eigenvalues, with non null imaginary part, crosses the imaginary axis (their real part becomes strictly positive). 
Even though the behaviour is qualitatively comparable, the presence of self- and/or cross-diffusion dramatically changes the value of the parameter of Hopf bifurcation.

%%%%%%%%%%%%%%%%%%%%%%%%%%%% Hopf %%%%%%%%%%%%%%%%%%%%%%%%%%%%%%%%%%%%%%

\subsection{Bifurcation of periodic orbits}

\subsubsection{Periodic orbits}

Let us denote by $\phi^t(X_0)$ the orbit, at time $t$, starting from $X_0 = (u_1(t=0,x),u_2(t=0,x))$ at time $t=0$, solution of the system (\ref{model}).
Periodic orbits are orbits $\bar{X}$ for which there exists a period $T>0$ such that:
\begin{equation}
\forall t, \phi^{T}(\bar{X}(t)) = \bar{X}(t).
\label{periodic_orbit}
\end{equation}
Computing a periodic orbit consists of finding a period $T$, and a corresponding limit cycle $\bar{X}$, solving equation (\ref{periodic_orbit}).
However, this system is not closed. 
First, $T$ is unknown a priori, and, even for a fixed $T$, there exists a continuum of solution, parametrized by a shift of time $\sigma$.
In order to close the problem, a phase condition has to be added: this is the basis of the shooting method \cite{Kuznetsov:2004}. 
The classical method consists of minimizing the distance $\rho$ between the periodic orbit $\bar{X}$ and the reference solution $\tilde{X}$, including a potential time shift $\sigma$:
\[
\rho(\sigma) = \int_0^T \|\tilde{X}(\tau + \sigma) - \bar{X}(\tau) \|^2 \mathrm{d}\tau.
\]
Considering a variation of $\sigma$ leads to the first order condition, cancelling the derivative:
\[
\rho'(\sigma) = 2\int_0^T (\tilde{X}(\tau + \sigma) - \bar{X}(\tau)).F(\tilde{X}(\tau)) \mathrm{d}\tau
\]
In our case, we consider, around $\tau=0$, with $\sigma = T$:
\[
\rho'(T) \approx (\phi^T(\bar{X}(0)) - \bar{X}(0).F(\tilde{X}(0)) ).
\]
The latter approximation is rather coarse. However, it is fast to evaluate, and works well in practice.
On our system, at least, it seems to increase the convergence of the solver.
Let us assume that we know a reference solution $\tilde{X} = (\tilde{u}_1,\tilde{u}_2)$. 
In practice, this will be the solution found using the former parameter value in a continuation loop, and the steady state before the Hopf bifurcation for the initial one. 
The periodic orbit problem may then be solved using a Newton-Krylov method on the system:
\begin{align}
\phi^{T}(\bar{X}(t)) - \bar{X}(t) = 0,\\
(\phi^T(\bar{X}(0)) - \bar{X}(0).F(\tilde{X}(0)) ) = 0.
\label{shooting}
\end{align}
For our tests, we have set up a tolerance of $5e^{-4}$ on the overall residual.

\subsubsection{Floquet theory and orbit stability}

Once the periodic orbit $\bar{X}$ is known, together with its period $T$, comes the question of its stability. 
The classical tool to do so is the concept of Poincaré map,
sending points from a section in the phase space (called the Poincaré section), to itself, after a period $T$. 
In particular, periodic orbits are fixed points of the Poincaré map.
Monodromy matrices $M$ (from the Greek "monos" and "dromos", "one road" matrix) are linear approximations of the Poincaré map, defined as:
\begin{equation}
M = \frac{\partial \phi^T(\bar{X}(0))}{\partial x}.
\label{monodromy}
\end{equation}
In practice, we will evaluate $M= (c_1,...c_{2N})$ by columns, by a finite difference approximation:
\[
c_i \approx \frac{\phi^T(\bar{X}(0) + h e_i) - \phi^T(\bar{X}(0)))}{h}.
\]
In our applications, we set $h = 0.001$.
By application of the Floquet theorem \cite{Floquet:1883}, 
we know that at each eigenvalue $\mu$ of the monodromy matrix, called a Floquet multipliers, corresponds a mode $X_{\mu}$, such that:
\begin{equation}
\dot{X}_{\mu} = \mu X_{\mu}.
\label{Floquet}
\end{equation}
Therefore, a mode is stable if and only if the absolute value of the Floquet multiplier $\mu$ is strictly smaller than $1$. 
Stability may be lost by 
a fold bifurcation ($\mu$ exit the unit circle at $1$),
a period doubling bifurcation ($\mu$ exits the unit circle at $-1$), leading to the creation of a cycle of double period,
or a Neimark Sacker bifurcation ($\mu$ exits the unit circle at $e^{i\theta}$), leading to the creation of a quasiperiodic cycle (2 periods coexist, but are not commensurate).

\subsubsection{Results}

Using a numerical continuation approach \cite{Kuznetsov:2004}, we study the orbits in each of the four cases, starting from the Hopf bifurcation point, 
where a limit cycle is created.
For each case, we display 
the variation of density of the first species $u_1(t)$, 
the corresponding variation of its Fourier transform in space $\frac{2}{N} |\mathcal{F}(u - \bar{u})|$, 
the evolution of energy (\ref{E}), 
the Fourier spectrum of the signal of the energy $\frac{2}{N} |\mathcal{F}(E - \bar{E})|$, 
the Floquet multipliers when the cycle becomes unstable,
and the periods with respect to the bifurcation parameter.

\begin{remark}
With this choice of convention in the plot of the Fourier transform, 
a signal of the form $E(t) = A\cos(2\pi f t)$ (resp. $u(x) = A\cos(k x)$) 
will appear as a peak on the Fourier domain, 
centered on $f=n$ (resp $k=n$), with height $A$, and therefore will be easily interpretable. 
Removing the mean values allows to define a refined color bar.
\end{remark} 

\begin{figure}[!htbp]\centering
%%%%%%%%%%%%%%%%%%%%%TIME DENSITY%%%%%%%%%%%%%%%%%%%%%%%%%%%%%%
\begin{subfigure}{0.24\textwidth}
\includegraphics[width=\linewidth]{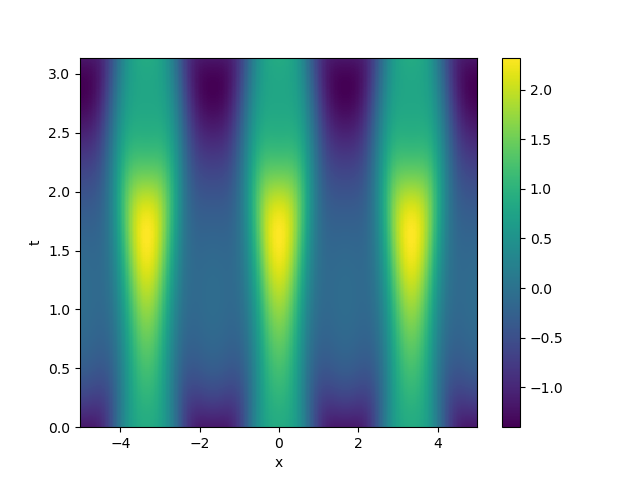}
\end{subfigure}
\begin{subfigure}{0.24\textwidth}
\includegraphics[width=\linewidth]{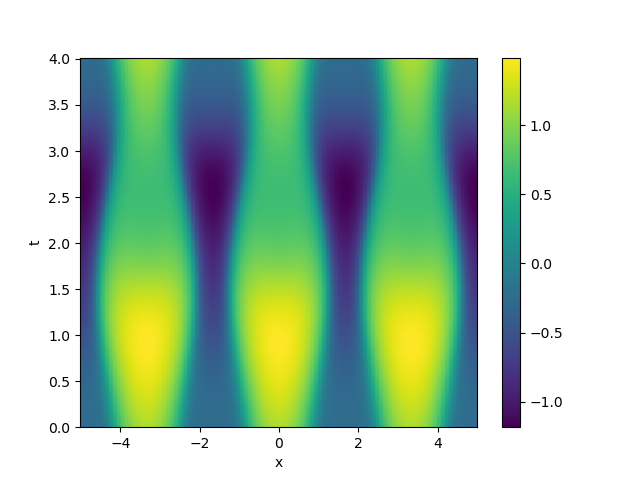}
\end{subfigure}
\begin{subfigure}{0.24\textwidth}
\includegraphics[width=\linewidth]{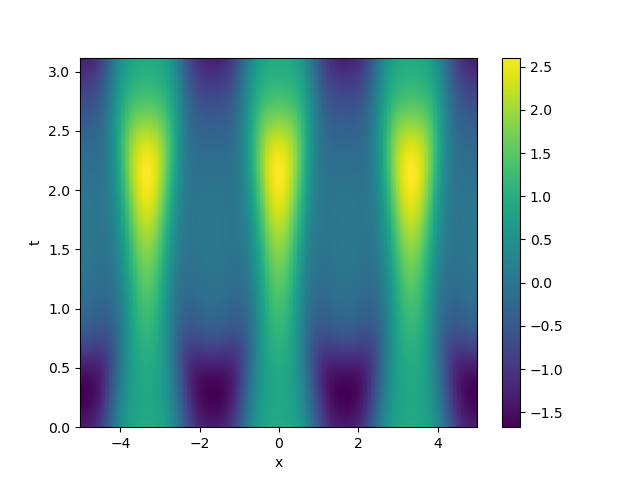}
\end{subfigure}
\begin{subfigure}{0.24\textwidth}
\includegraphics[width=\linewidth]{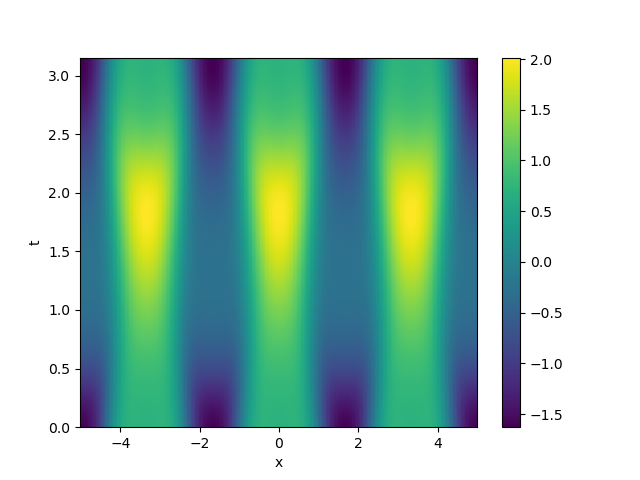}
\end{subfigure}\\
%%%%%%%%%%%%%%%%%%%%%TIME FOURIER%%%%%%%%%%%%%%%%%%%%%%%%%%%%%%
\begin{subfigure}{0.24\textwidth}
\includegraphics[width=\linewidth]{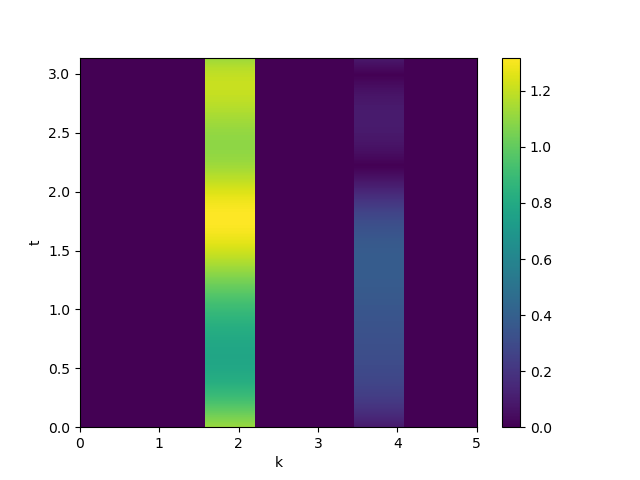}
\end{subfigure}
\begin{subfigure}{0.24\textwidth}
\includegraphics[width=\linewidth]{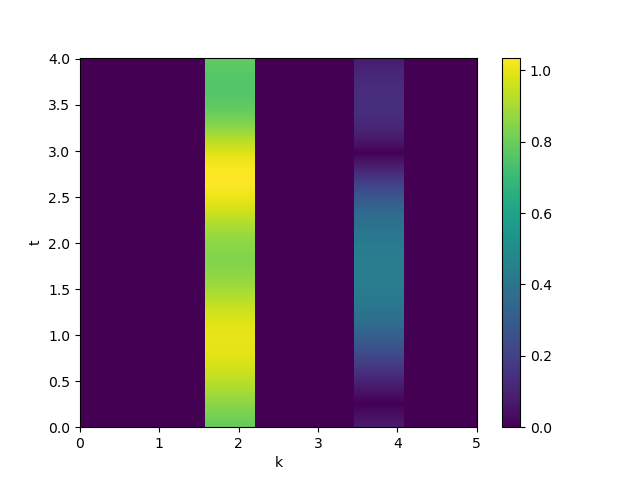}
\end{subfigure}
\begin{subfigure}{0.24\textwidth}
\includegraphics[width=\linewidth]{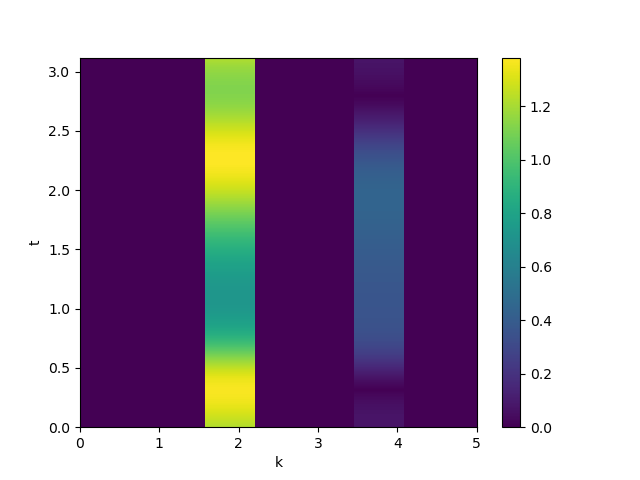}
\end{subfigure}
\begin{subfigure}{0.24\textwidth}
\includegraphics[width=\linewidth]{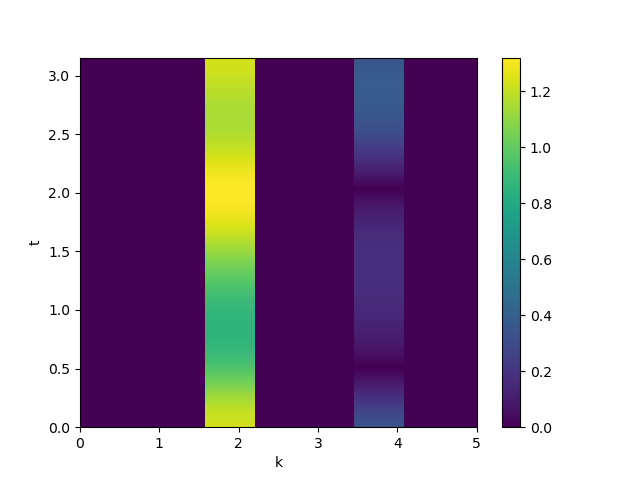}
\end{subfigure}\\
%%%%%%%%%%%%%%%%%%%%%%%%%Energy%%%%%%%%%%%%%%%%%%%%%%%%%%%%%
\begin{subfigure}{0.24\textwidth}
\includegraphics[width=\linewidth]{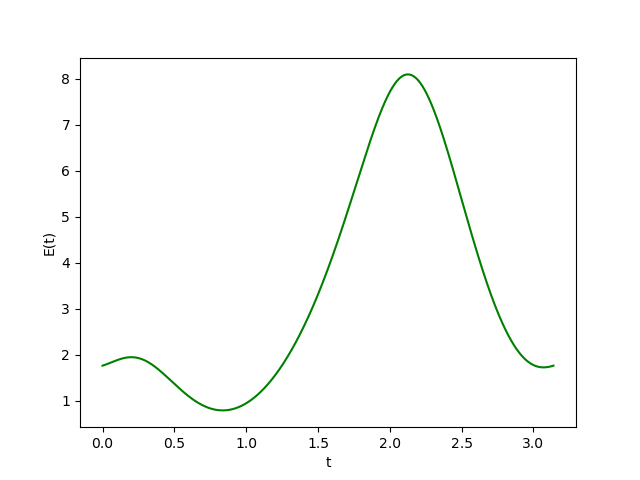}
\end{subfigure}
\begin{subfigure}{0.24\textwidth}
\includegraphics[width=\linewidth]{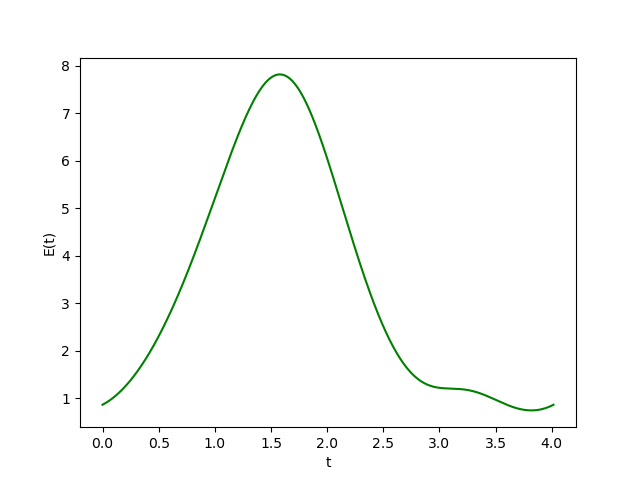}
\end{subfigure}
\begin{subfigure}{0.24\textwidth}
\includegraphics[width=\linewidth]{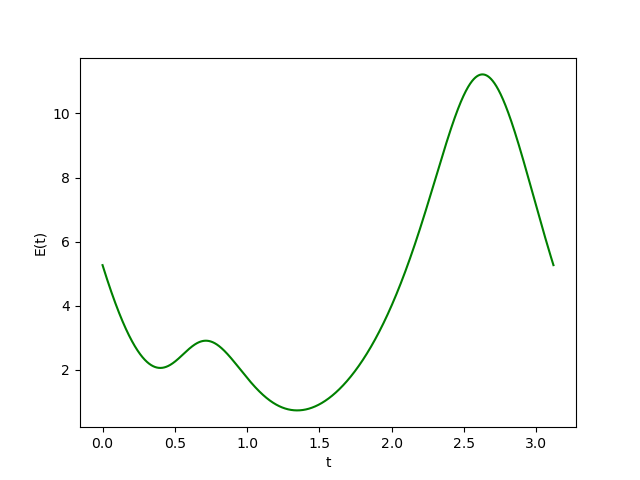}
\end{subfigure}
\begin{subfigure}{0.24\textwidth}
\includegraphics[width=\linewidth]{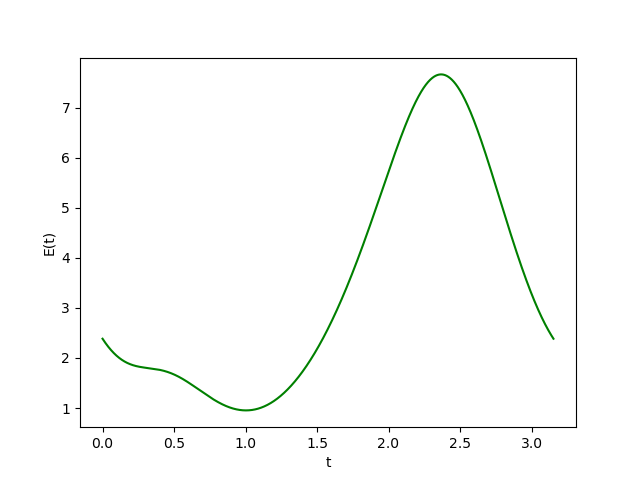}
\end{subfigure}\\
%%%%%%%%%%%%%%%%%%%%%%%%%Spectrum%%%%%%%%%%%%%%%%%%%%%%%%%%%%%
\begin{subfigure}{0.24\textwidth}
\includegraphics[width=\linewidth]{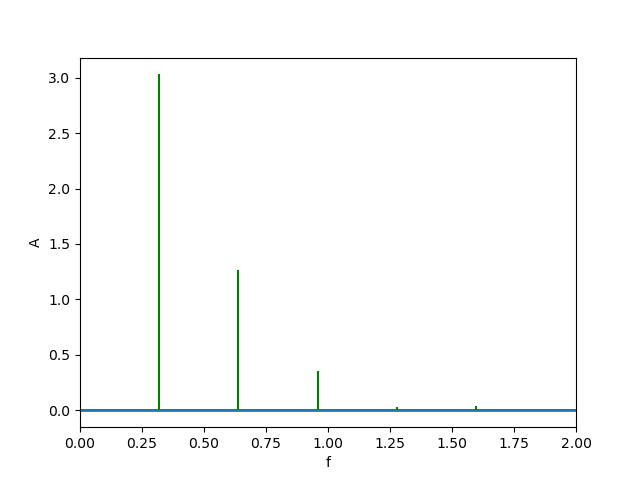}
\end{subfigure}
\begin{subfigure}{0.24\textwidth}
\includegraphics[width=\linewidth]{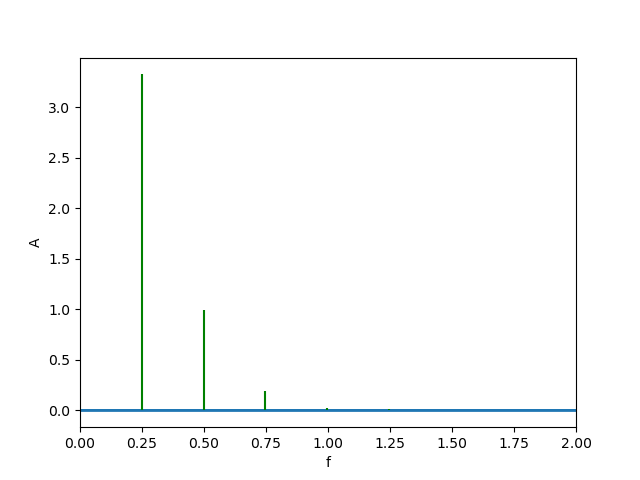}
\end{subfigure}
\begin{subfigure}{0.24\textwidth}
\includegraphics[width=\linewidth]{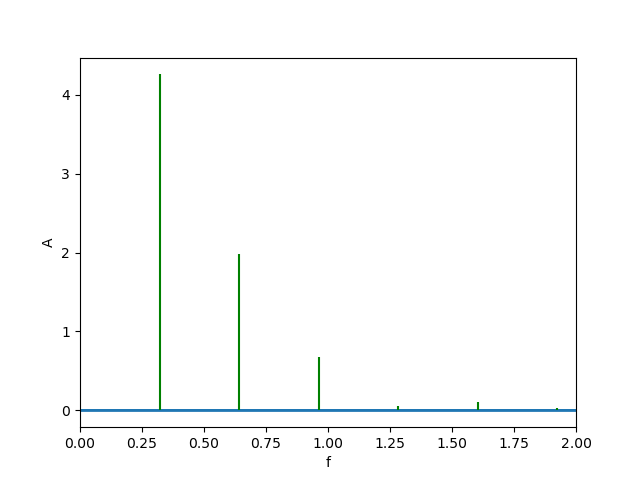}
\end{subfigure}
\begin{subfigure}{0.24\textwidth}
\includegraphics[width=\linewidth]{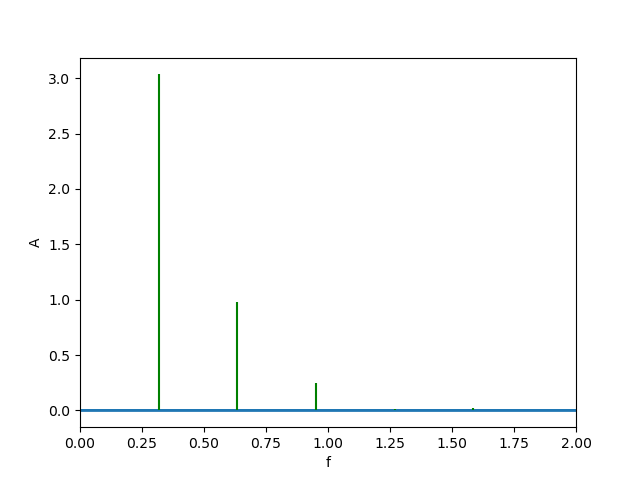}
\end{subfigure}\\
%%%%%%%%%%%%%%%%%%%%%Period%%%%%%%%%%%%%%%%%%%%%%%%%%%%%%
\begin{subfigure}{0.24\textwidth}
\includegraphics[width=\linewidth]{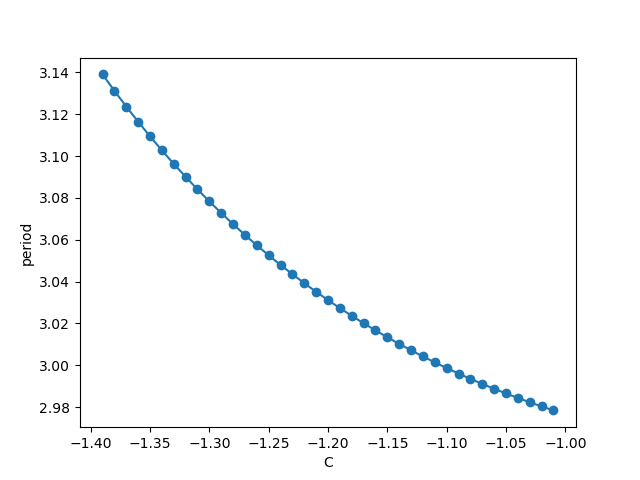}
\end{subfigure}
\begin{subfigure}{0.24\textwidth}
\includegraphics[width=\linewidth]{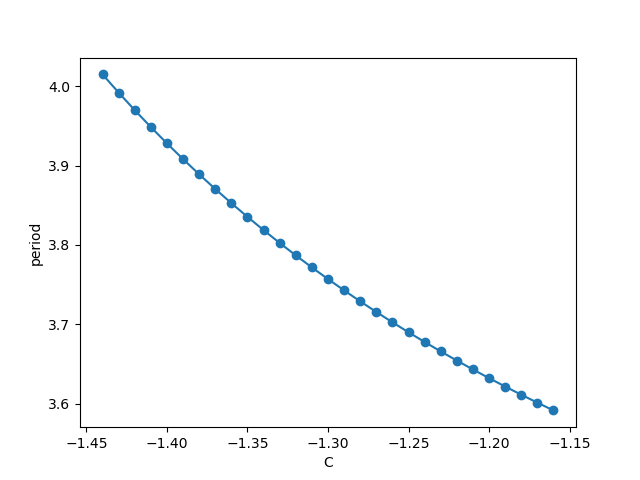}
\end{subfigure}
\begin{subfigure}{0.24\textwidth}
\includegraphics[width=\linewidth]{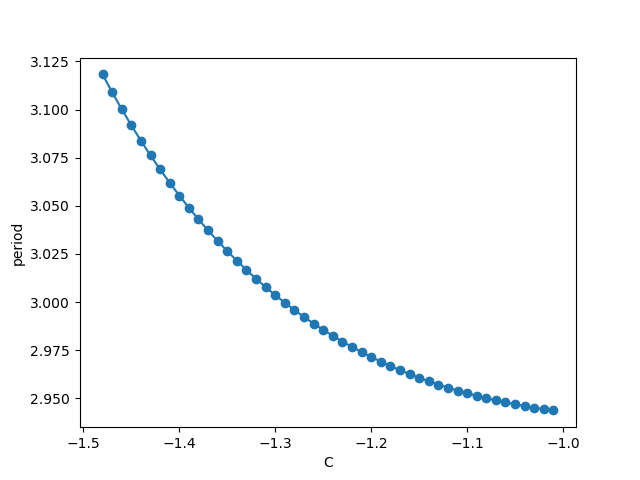}
\end{subfigure}
\begin{subfigure}{0.24\textwidth}
\includegraphics[width=\linewidth]{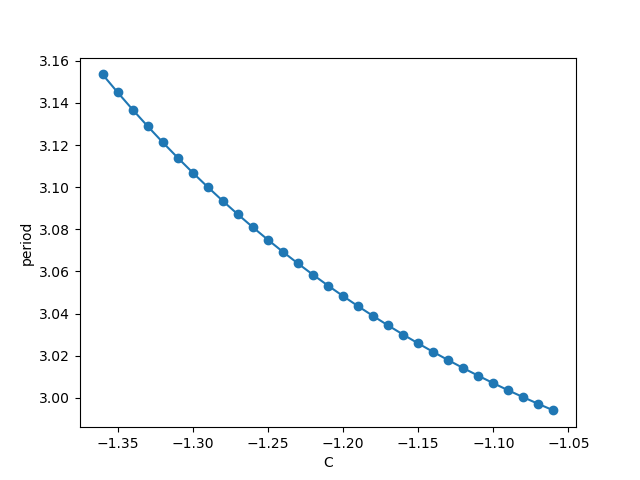}
\end{subfigure}\\
%%%%%%%%%%%%%%%%%%%%%Floquet%%%%%%%%%%%%%%%%%%%%%%%%%%%%%%
\begin{subfigure}{0.24\textwidth}
\includegraphics[width=\linewidth]{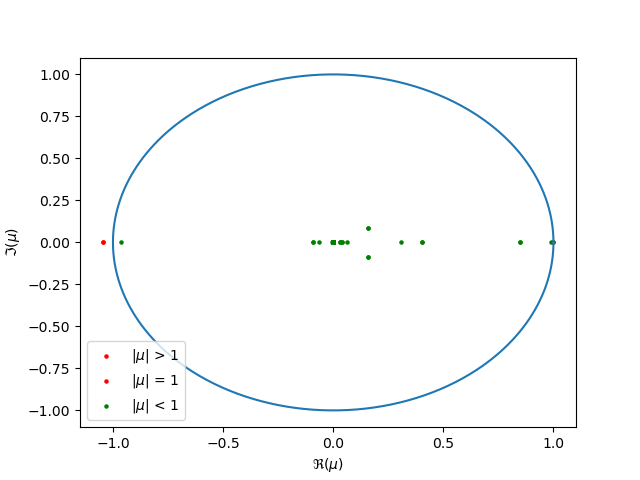}
\end{subfigure}
\begin{subfigure}{0.24\textwidth}
\includegraphics[width=\linewidth]{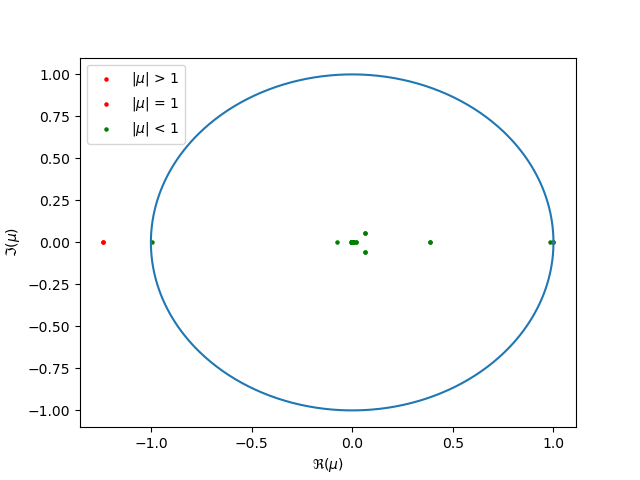}
\end{subfigure}
\begin{subfigure}{0.24\textwidth}
\includegraphics[width=\linewidth]{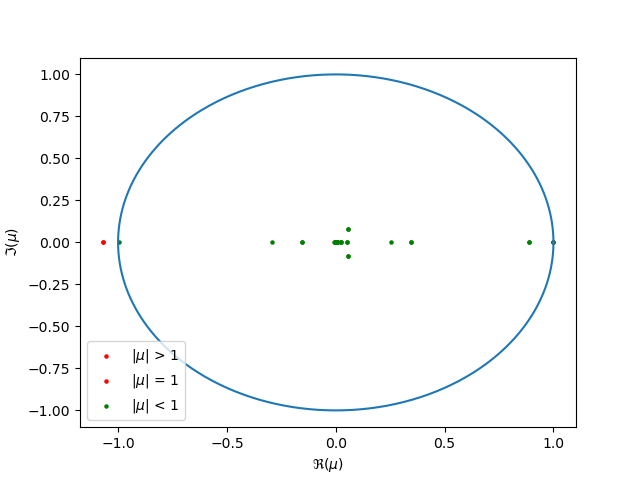}
\end{subfigure}
\begin{subfigure}{0.24\textwidth}
\includegraphics[width=\linewidth]{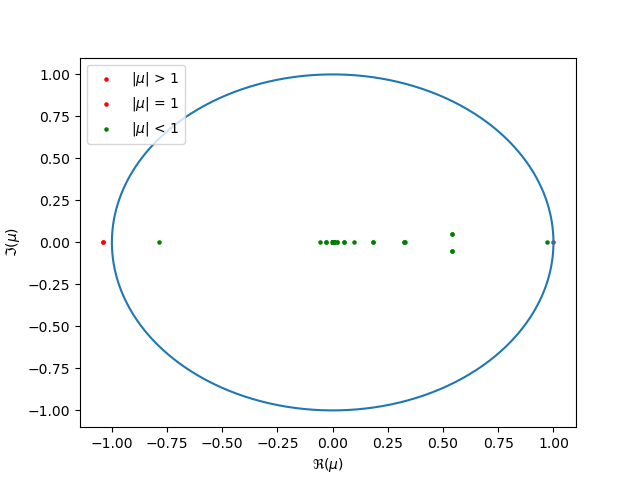}
\end{subfigure}\\
\caption{Bifurcation of periodic orbits, starting from the Hopf bifurcation (Figure \ref{Hopf}).
Periodic orbits are computed using the shooting method (\ref{shooting}). 
Density and its Fourier transform (in space) are computed with respect to time, and one full cycle is displayed on lines 1 and 2.
Energy (\ref{E}) and its Fourier transform (in time), on one full cycle, are displayed on lines 3 and 4.
Period with respect to control parameter C is displayed on line 5.
Floquet multipliers (\ref{Floquet}) are displayed on line 6, at  $C = -1.39, -1.43, -1.48, -1.36$.
Eventually, all periodic orbits lose their stability via a period doubling bifurcation 
($\mu$ exits the unit circle via $-1$).}
\label{1cycle}
\end{figure}

On Figure \ref{1cycle}, we consider a variation of $C$ with step $\delta C = -0.01$.
We observe that,
in each of the four cases, passed the value of the Hopf bifurcation, a periodic orbit appears, 
in perfect agreement with the value computed in the previous section (see Figure \ref{Hopf}).
Starting from the unstable steady state, and after an initial simulation (to be sure to achieve convergence), 
we use this periodic orbit as inital condition of our continuation strategy.
As initial guess, we set $T = 3$, according to observations.
The Turing patterns, once steady, become unstable, and start to oscillate in time.
For each of the four cases, a period doubling bifurcation is observed, when the value of parameter $C$ is decreased,
as a Floquet multiplier, in red, crosses the unit circle at $-1$.
We notice a peak around $k=2$ for the density in Fourier domain, leading to a spatial frequency of $f = \frac{2 \pi}{k} \approx \pi$, 
in agreement with the distribution of density in space.

\begin{figure}[!htbp]\centering
%%%%%%%%%%%%%%%%%%%%%TIME DENSITY%%%%%%%%%%%%%%%%%%%%%%%%%%%%%%
\begin{subfigure}{0.24\textwidth}
\includegraphics[width=\linewidth]{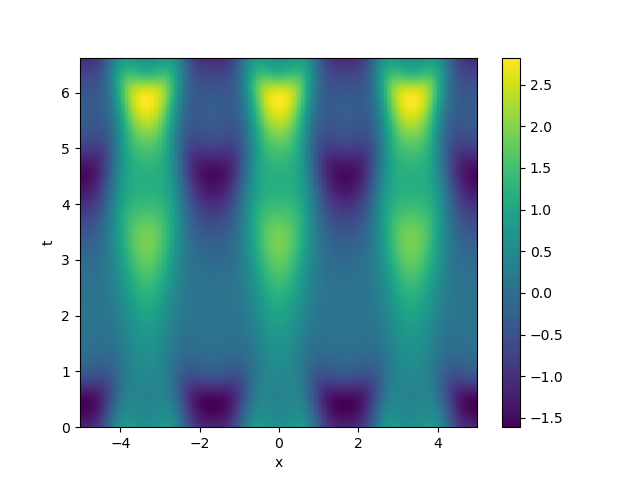}
\end{subfigure}
\begin{subfigure}{0.24\textwidth}
\includegraphics[width=\linewidth]{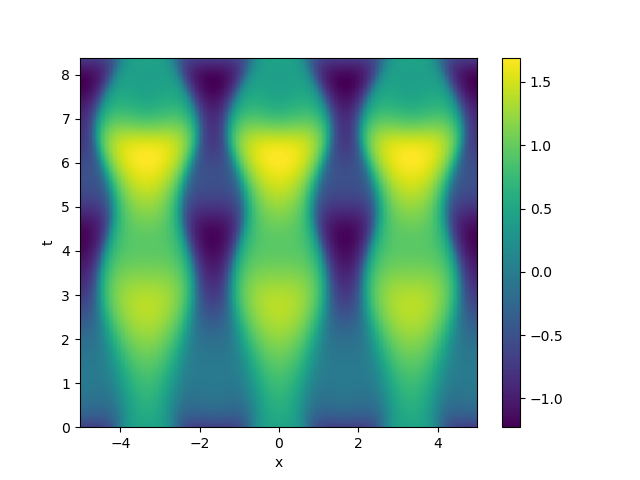}
\end{subfigure}
\begin{subfigure}{0.24\textwidth}
\includegraphics[width=\linewidth]{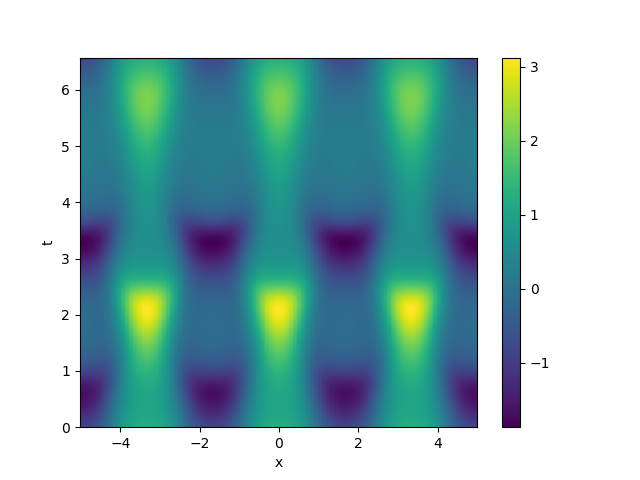}
\end{subfigure}
\begin{subfigure}{0.24\textwidth}
\includegraphics[width=\linewidth]{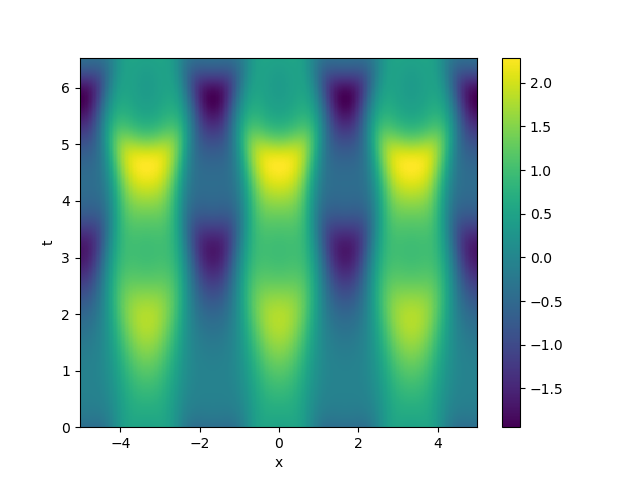}
\end{subfigure}\\
%%%%%%%%%%%%%%%%%%%%%TIME FOURIER%%%%%%%%%%%%%%%%%%%%%%%%%%%%%%
\begin{subfigure}{0.24\textwidth}
\includegraphics[width=\linewidth]{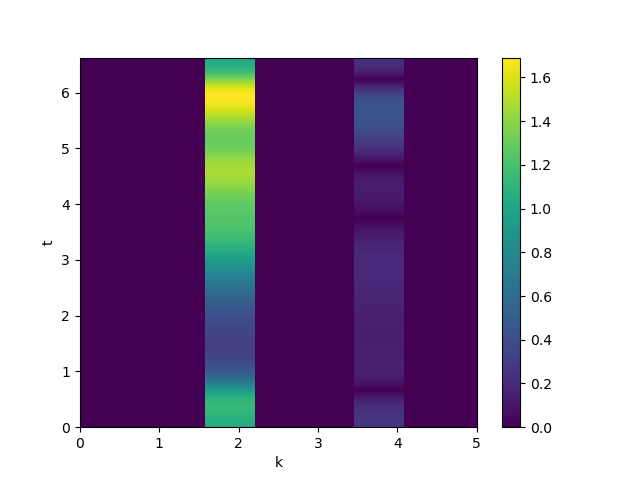}
\end{subfigure}
\begin{subfigure}{0.24\textwidth}
\includegraphics[width=\linewidth]{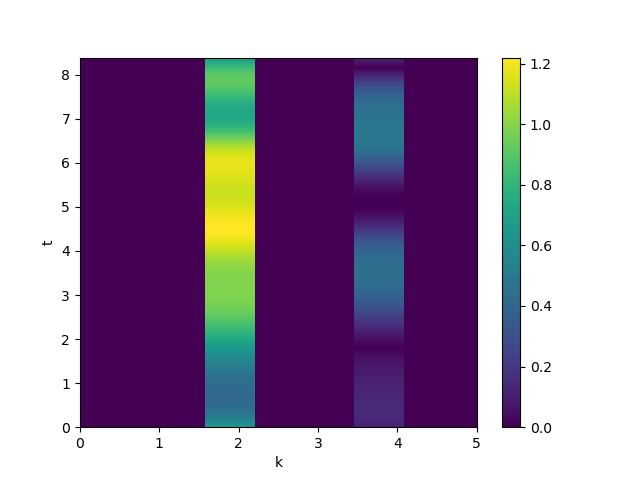}
\end{subfigure}
\begin{subfigure}{0.24\textwidth}
\includegraphics[width=\linewidth]{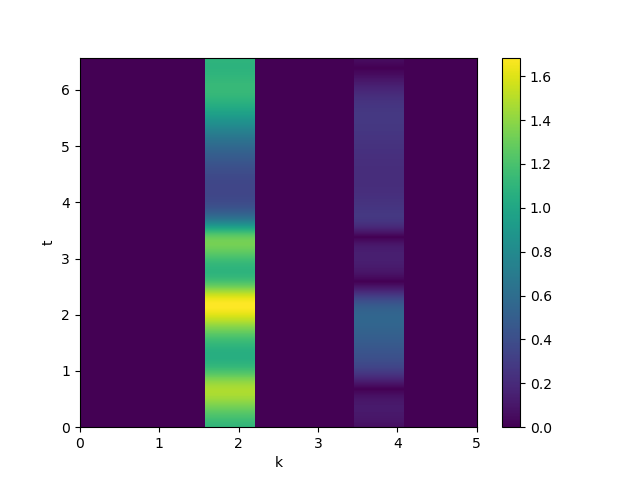}
\end{subfigure}
\begin{subfigure}{0.24\textwidth}
\includegraphics[width=\linewidth]{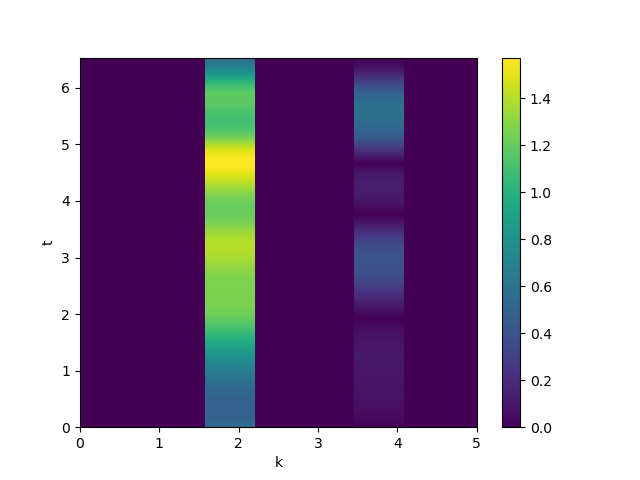}
\end{subfigure}\\
%%%%%%%%%%%%%%%%%%%%%%%%%Energy%%%%%%%%%%%%%%%%%%%%%%%%%%%%%
\begin{subfigure}{0.24\textwidth}
\includegraphics[width=\linewidth]{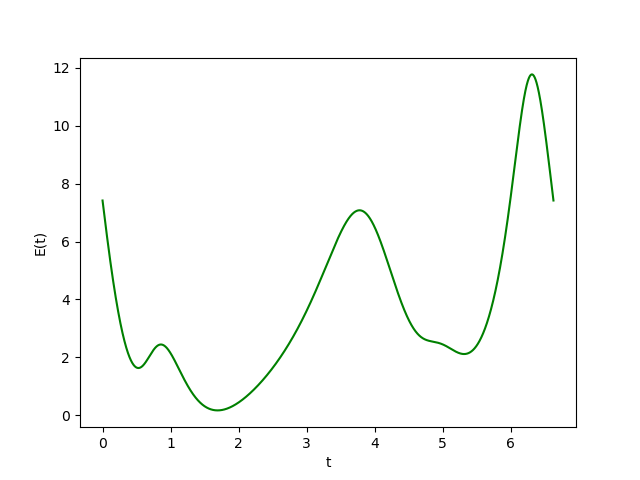}
\end{subfigure}
\begin{subfigure}{0.24\textwidth}
\includegraphics[width=\linewidth]{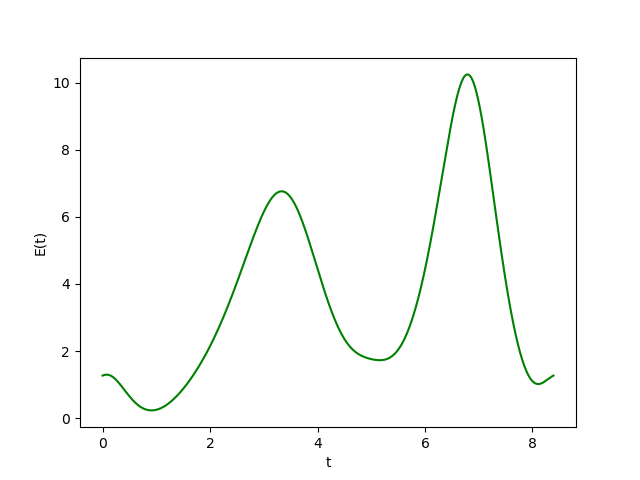}
\end{subfigure}
\begin{subfigure}{0.24\textwidth}
\includegraphics[width=\linewidth]{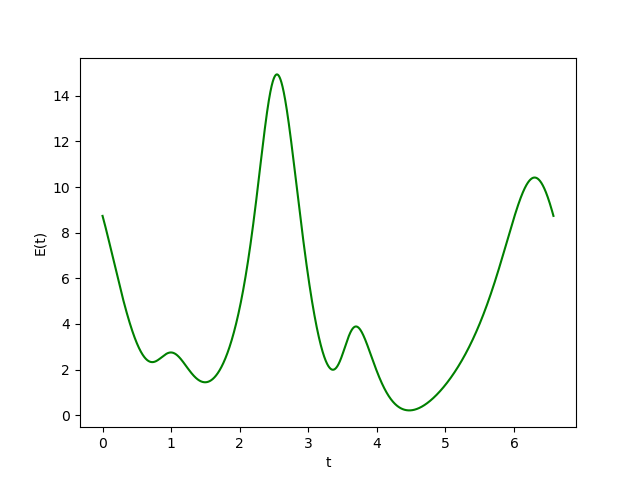}
\end{subfigure}
\begin{subfigure}{0.24\textwidth}
\includegraphics[width=\linewidth]{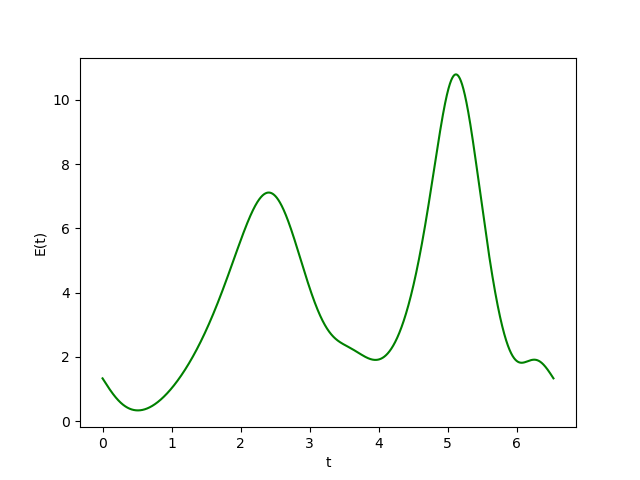}
\end{subfigure}\\
%%%%%%%%%%%%%%%%%%%%%%%%%Spectrum%%%%%%%%%%%%%%%%%%%%%%%%%%%%%
\begin{subfigure}{0.24\textwidth}
\includegraphics[width=\linewidth]{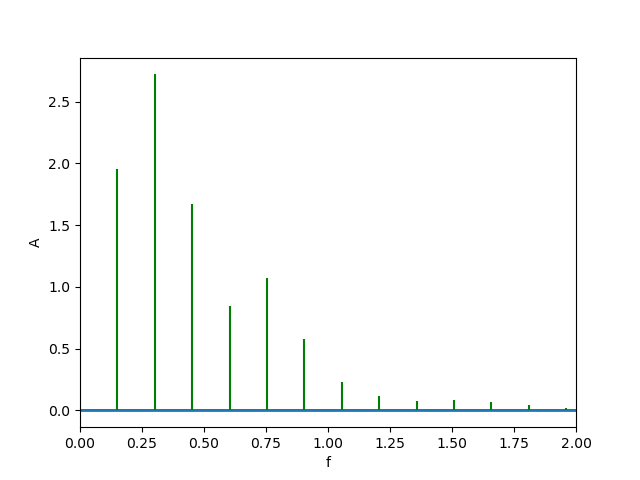}
\end{subfigure}
\begin{subfigure}{0.24\textwidth}
\includegraphics[width=\linewidth]{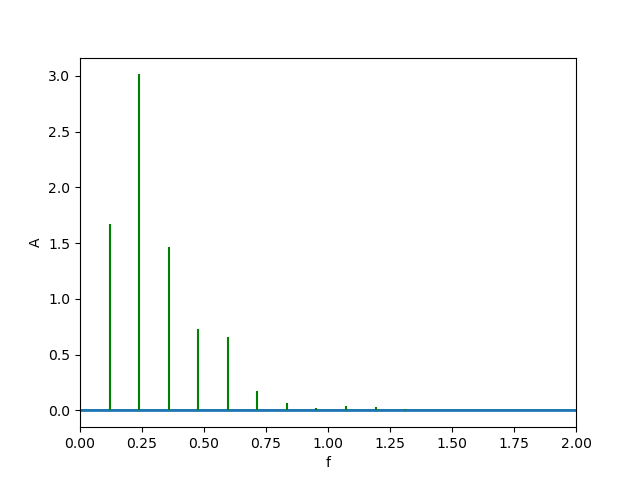}
\end{subfigure}
\begin{subfigure}{0.24\textwidth}
\includegraphics[width=\linewidth]{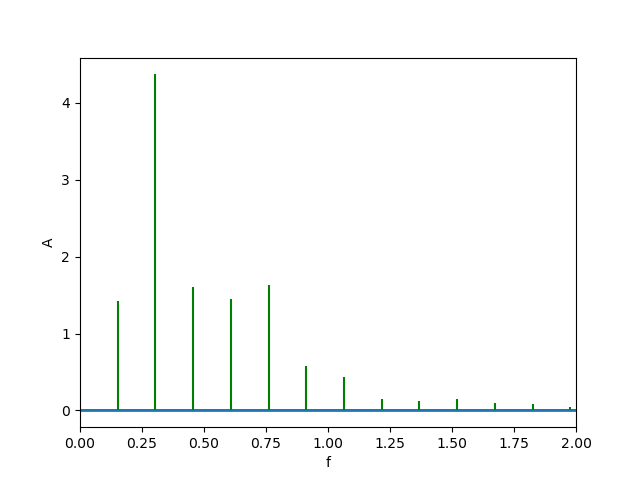}
\end{subfigure}
\begin{subfigure}{0.24\textwidth}
\includegraphics[width=\linewidth]{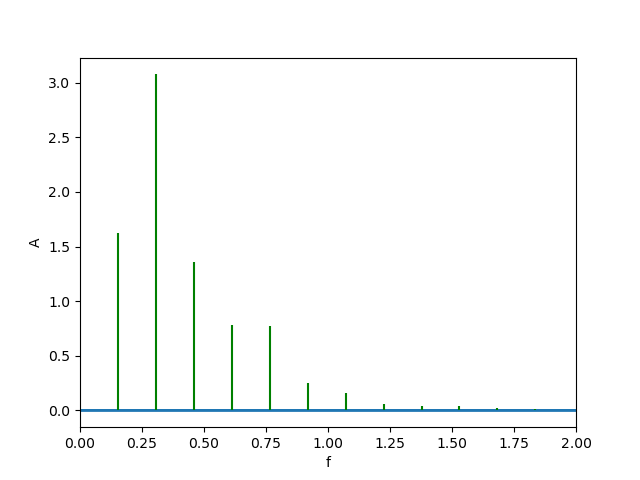}
\end{subfigure}\\
%%%%%%%%%%%%%%%%%%%%%Periods%%%%%%%%%%%%%%%%%%%%%%%%%%%%%%
\begin{subfigure}{0.24\textwidth}
\includegraphics[width=\linewidth]{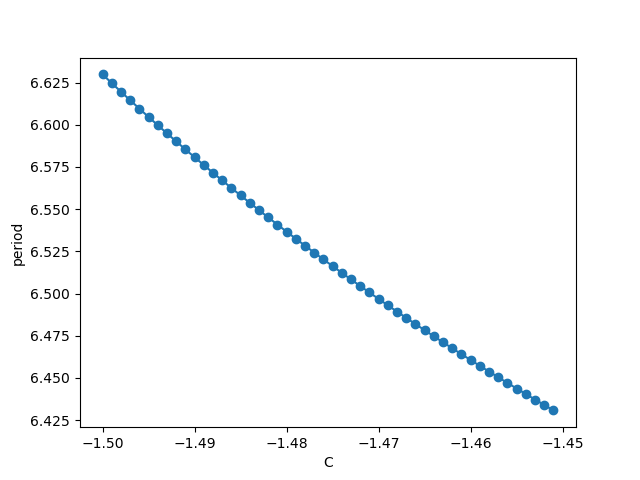}
\end{subfigure}
\begin{subfigure}{0.24\textwidth}
\includegraphics[width=\linewidth]{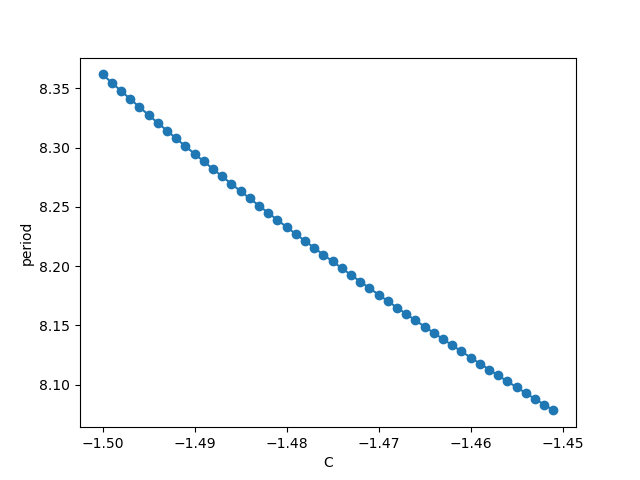}
\end{subfigure}
\begin{subfigure}{0.24\textwidth}
\includegraphics[width=\linewidth]{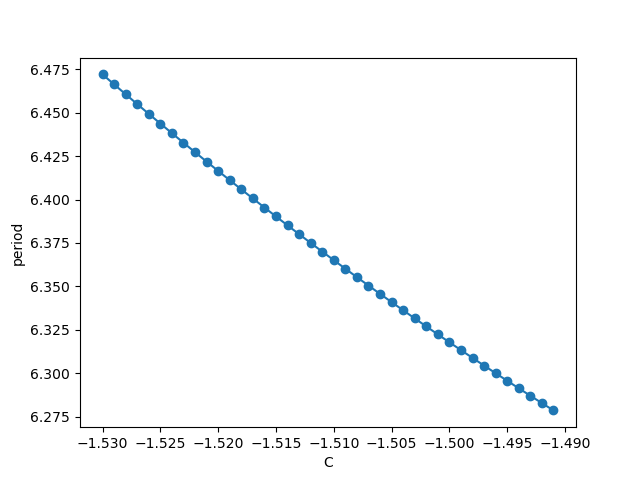}
\end{subfigure}
\begin{subfigure}{0.24\textwidth}
\includegraphics[width=\linewidth]{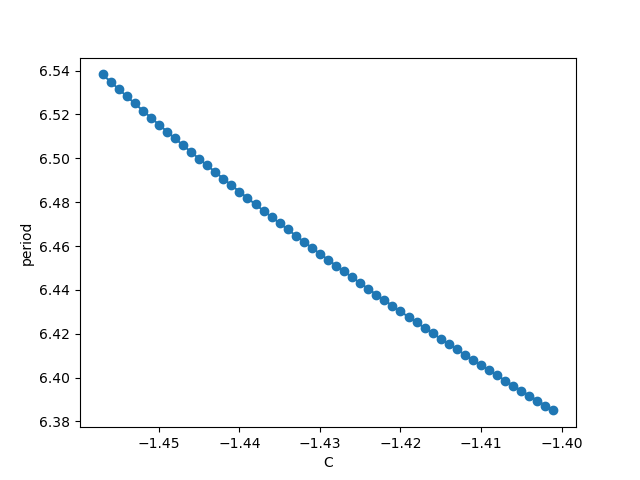}
\end{subfigure}\\
%%%%%%%%%%%%%%%%%%%%%Floquet%%%%%%%%%%%%%%%%%%%%%%%%%%%%%%
\begin{subfigure}{0.24\textwidth}
\includegraphics[width=\linewidth]{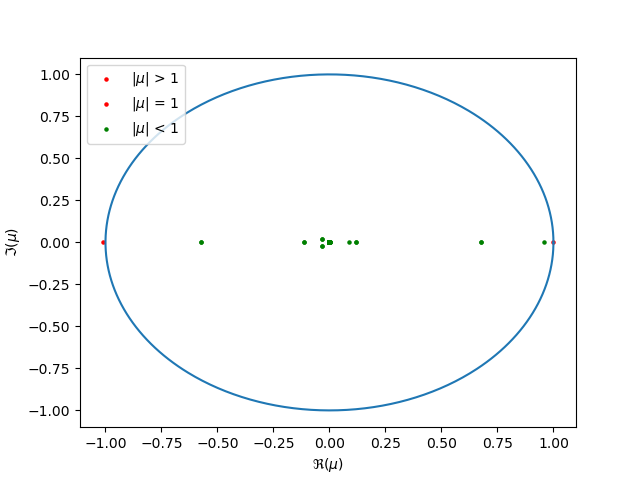}
\end{subfigure}
\begin{subfigure}{0.24\textwidth}
\includegraphics[width=\linewidth]{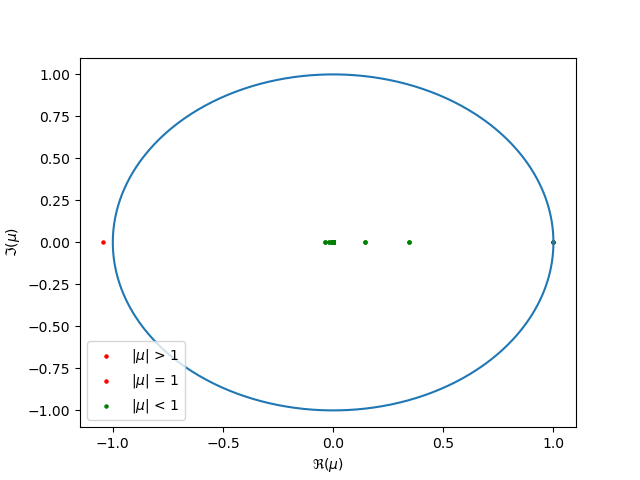}
\end{subfigure}
\begin{subfigure}{0.24\textwidth}
\includegraphics[width=\linewidth]{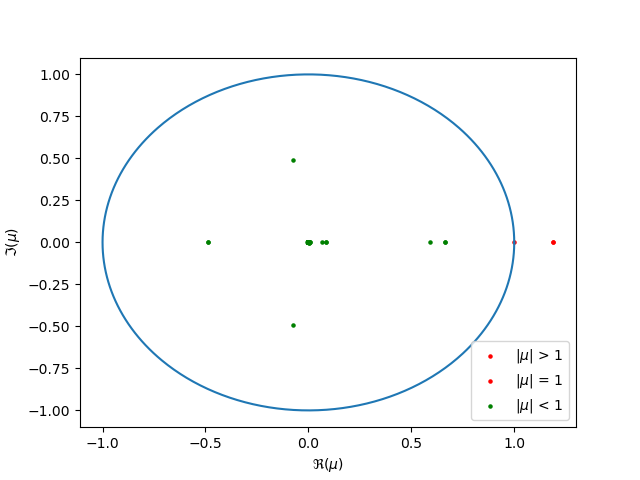}
\end{subfigure}
\begin{subfigure}{0.24\textwidth}
\includegraphics[width=\linewidth]{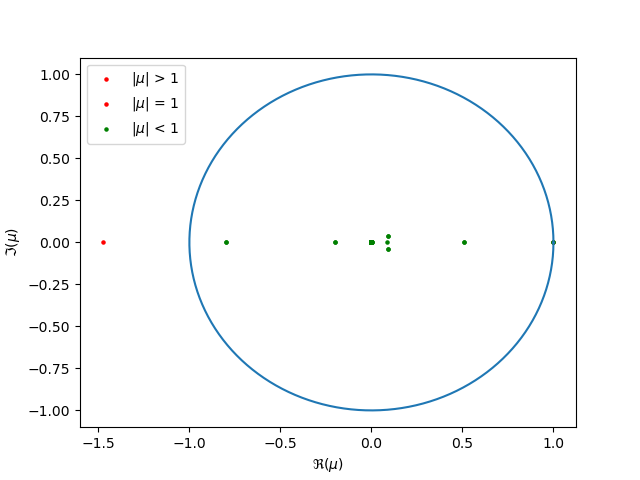}
\end{subfigure}\\
\caption{Bifurcation of periodic orbits, starting from the first period doubling bifurcation (Figure \ref{1cycle}).
Floquet multipliers are displayed on line 6, at  $C = -1.471, -1.482, -1.599, -1.486$.
Eventually, all periodic lose stability, either by period doubling ($\mu$ exits the unit circle via $-1$) or by a fold ($\mu$ exits via $1$).}
\label{2cycle}
\end{figure}

On Figure \ref{2cycle}, we start from the unstable periodic orbit, and we simulate the orbit until the periodic doubling is observed.
We use this orbit as an initial condition for the shooting method with an initial guess of $T$ as twice the value of the previous limit cycle.
We consider a smaller step of $\delta C = -0.001$, as the previous step was not refined enough to ensure convergence between consecutive values of birffurcation parameter.
We continue the orbits until they become unstable.
Interestingly, the scenario changes according to the case.
For linear diffusion, self diffusion on $u_1$, and cross-diffusion, a second period doubling is observed.
For the self diffusion on $u_2$, a fold bifurcation appears.

\begin{remark}
In order to check the convergence of numerical results, we have tried $N=300$ and $N=500$ in the former computations.
We have observed, qualitatively, the same results.
\end{remark}

%%%%%%%%%%%%%%%%%%%%%%% 1-cycles %%%%%%%%%%%%%%%%%%%%%%%%%%%%%%%%%%%%%%%%

%%%%%%%%%%%%%%%%%%%%%%% 2 cycles %%%%%%%%%%%%%%%%%%%%%%%%%%%%%%%%%%%%%%%%

\section{Chaos and strange attractors}

\subsection{Phase space}

In the original BVAM article \cite{BVAM:2012}, authors consider attractors in phase space of the form:
\begin{equation}
A(t) = \left\{(u_1(x=0,t),u_2(x=0,t)| t \in [\tau,T] \right\}.
\label{attractor_local}
\end{equation} 
In this work, we also monitor a phase space in terms of energy, by defining:
\begin{equation}
A_E(t) = \left\{E(t),\frac{d}{dt}E(t)| t \in [\tau,T] \right\}.
\label{attractor_global}
\end{equation}
This attractor is more global, in the sense that, contrarily to the previous one, it is not attached to a specific space coordinate ($x=0$), 
but only to a time $t$.
In addition, we color the curves with respect to a normalized time (sending $[0,T]$ onto $[0,1]$);
a blue color indicates a small time, while a yellow one a large time.

\subsection{Road to chaos}

\begin{figure}[!htbp]\centering
\begin{subfigure}{0.24\textwidth}
\includegraphics[width=\linewidth]{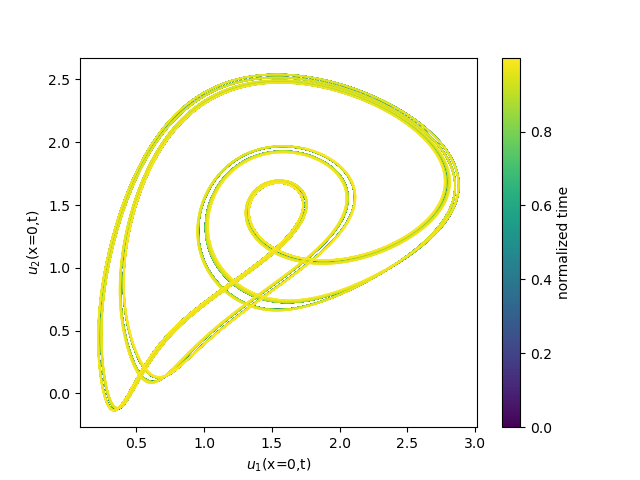}
\end{subfigure}
\begin{subfigure}{0.24\textwidth}
\includegraphics[width=\linewidth]{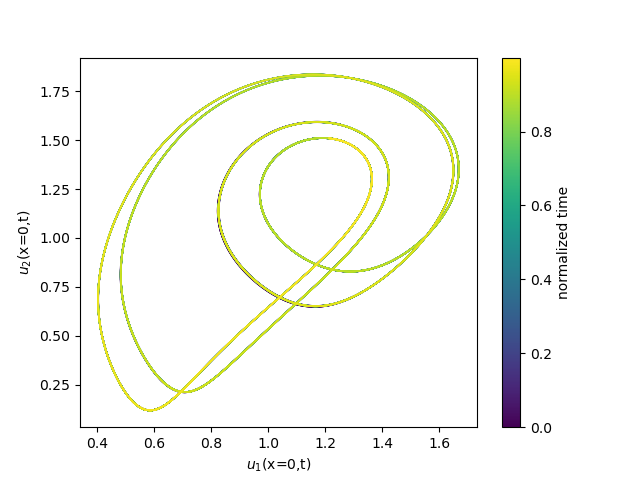}
\end{subfigure}
\begin{subfigure}{0.24\textwidth}
\includegraphics[width=\linewidth]{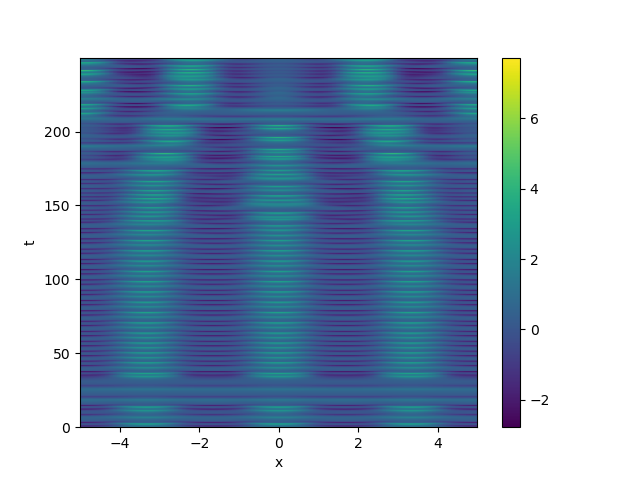}
\end{subfigure}
\begin{subfigure}{0.24\textwidth}
\includegraphics[width=\linewidth]{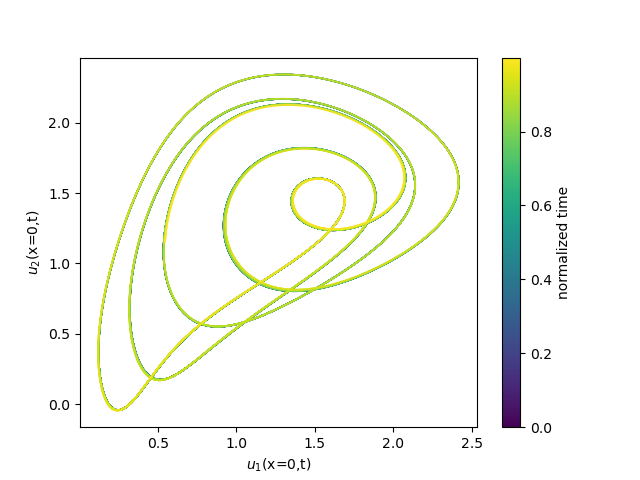}
\end{subfigure}\\
\caption{Road to chaos, from the instability of the periodic orbit (Figure \ref{2cycle}). 
$C = -1.531, -1.49, -1.63, -1.52$.
A period doubling cascade appears in three cases: linear diffusion, self-diffusion on $u_1$, and cross-diffusion, clearly visible in phase space, where the closed orbits double.}
\label{road}
\end{figure}

On Figure \ref{2cycle}, we observed that, for each of the four regimes of diffusion, the periodic orbit becomes unstable when the absolute value of the control parameter $C$ is high enough.
Eventually, for each case, the periodic orbit vanish, and the dynamics becomes chaotic.
However, the road to chaos may be different according to the diffusion regime, as illustrated on Figure \ref{road}.
For the cases of linear diffusion, self-diffusion on the first species, and cross-diffusion, a period doubling cascade to chaos is suspected;
indeed, several period doubling are catched by the periodic orbit solver, and are visible in phase space, as a closed curve, with multiple loops, with a yellow color (indicating that orbits are revisited, it is not a transient signal).
Let us note that the period doubling cascade observed in the linear case differs from the quasiperiodicity route orignally observed by BVAM. 
This could be due to the different boundary conditions under consideration: homogenous Dirichlet in the original case, as opposed to periodic condition in this article.
For the case of self-diffusion on the second species, a fold bifurcation is observed.
Observing the evolution of the density, it seems that the periodic orbit, consisting of three oscillating stripes, evolves towards a five stripes configuration, and, doing so, becomes chaotic.

%%%%%%%%%%%%%%%%%%%%%%% transition %%%%%%%%%%%%%%%%%%%%%%%%%%%%%%%%%%%%%%%%

\subsection{Strange attractors}

After a transient period $\tau$, the limit set $A$ can take various forms, depending on the value of parameter $C$.
Before the Hopf bifurcarion, $A$ is reduced to a point. 
Passed the bifurcation, a periodic orbit appears, then $A$ is a circle. 
Passed the period doubling bifurcation, $A$ is a connected sum of two embedded circles. 
In the chaotic regime, $A$ fills a volume of the phase space, forming a strange attractor.

\begin{figure}[!htbp]\centering
%%%%%%%%%%%%%%%%%%%%%%%%%PHASE%%%%%%%%%%%%%%%%%%%%%%%%%%%%%
\begin{subfigure}{0.3\textwidth}
\includegraphics[width=\linewidth]{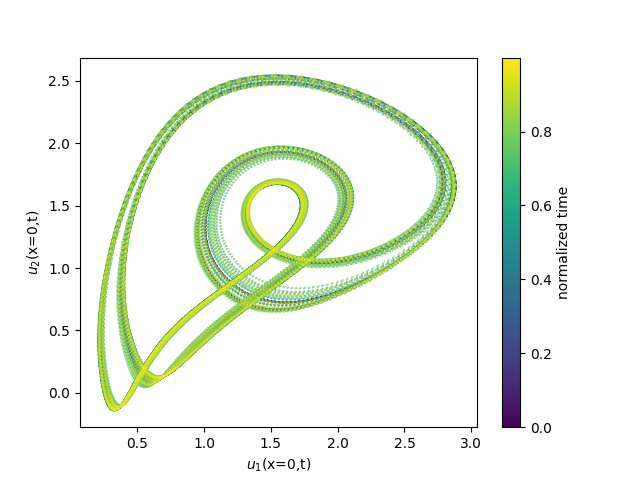}
\end{subfigure}
\begin{subfigure}{0.3\textwidth}
\includegraphics[width=\linewidth]{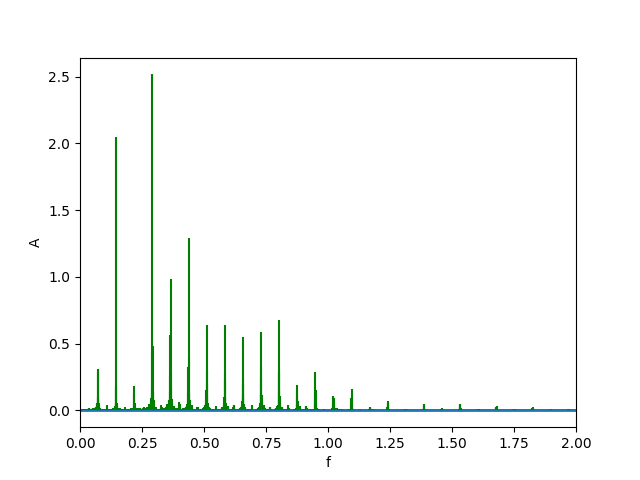}
\end{subfigure}
\begin{subfigure}{0.3\textwidth}
\includegraphics[width=\linewidth]{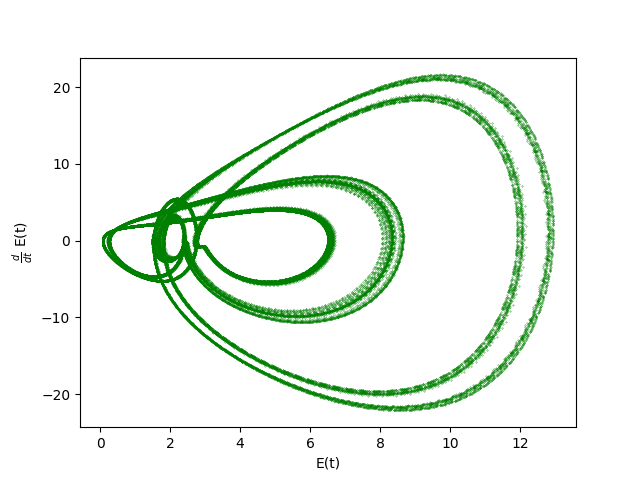}
\end{subfigure}
%%%%%%%%%%%%%%%%%%%%%%%%%PHASE%%%%%%%%%%%%%%%%%%%%%%%%%%%%%
\begin{subfigure}{0.3\textwidth}
\includegraphics[width=\linewidth]{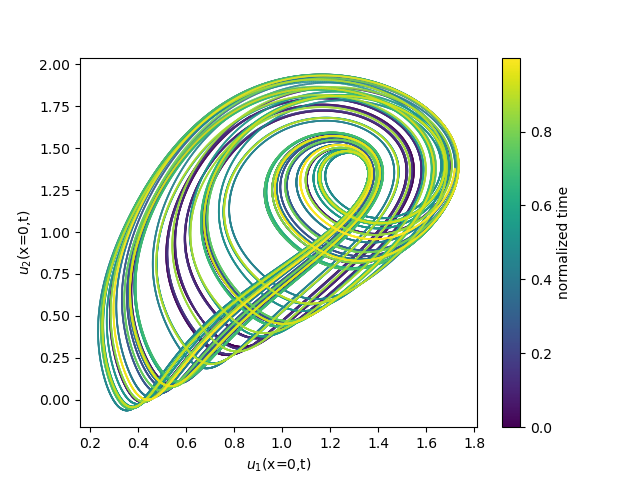}
\end{subfigure}
\begin{subfigure}{0.3\textwidth}
\includegraphics[width=\linewidth]{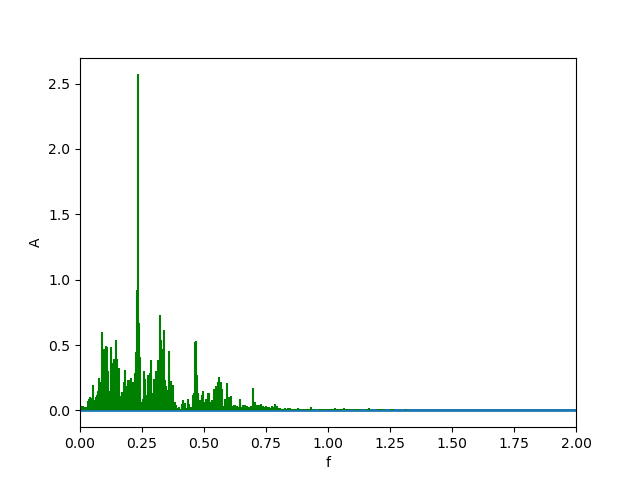}
\end{subfigure}
\begin{subfigure}{0.3\textwidth}
\includegraphics[width=\linewidth]{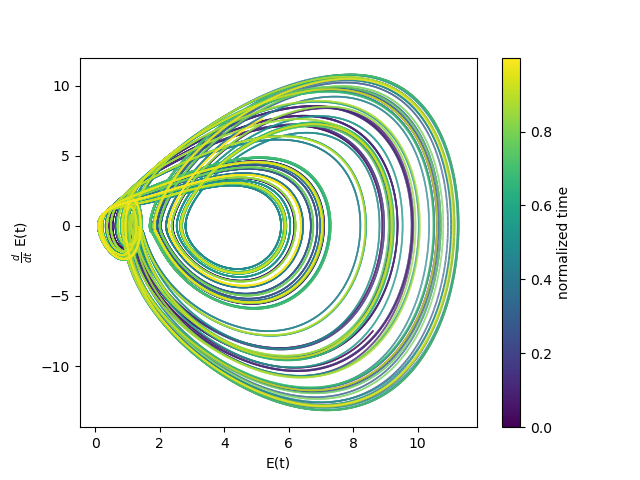}
\end{subfigure}
%%%%%%%%%%%%%%%%%%%%%%%%%PHASE%%%%%%%%%%%%%%%%%%%%%%%%%%%%%
\begin{subfigure}{0.3\textwidth}
\includegraphics[width=\linewidth]{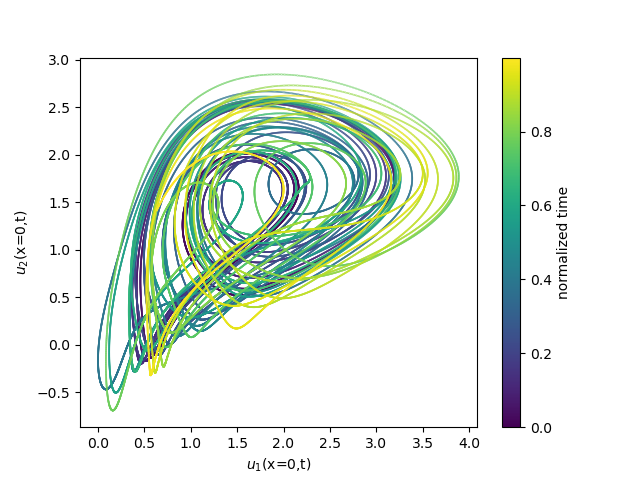}
\end{subfigure}
\begin{subfigure}{0.3\textwidth}
\includegraphics[width=\linewidth]{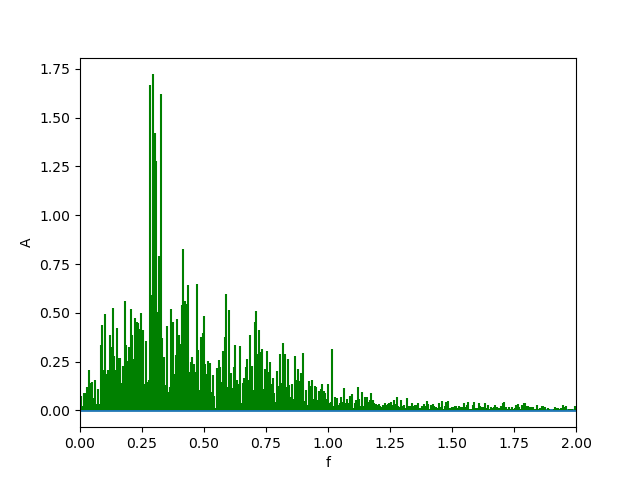}
\end{subfigure}
\begin{subfigure}{0.3\textwidth}
\includegraphics[width=\linewidth]{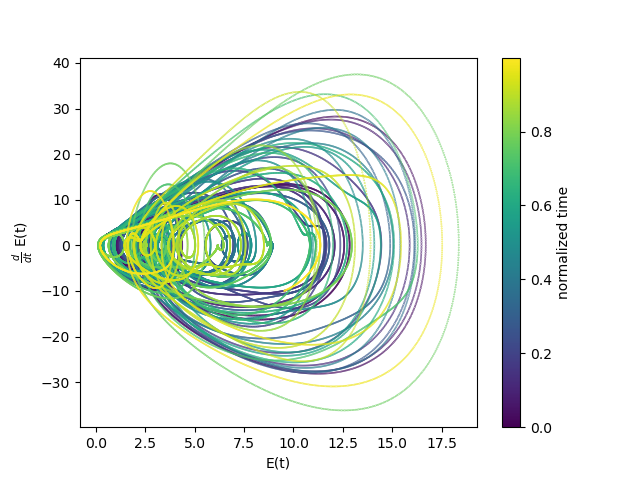}
\end{subfigure}
%%%%%%%%%%%%%%%%%%%%%%%%%PHASE%%%%%%%%%%%%%%%%%%%%%%%%%%%%%
\begin{subfigure}{0.3\textwidth}
\includegraphics[width=\linewidth]{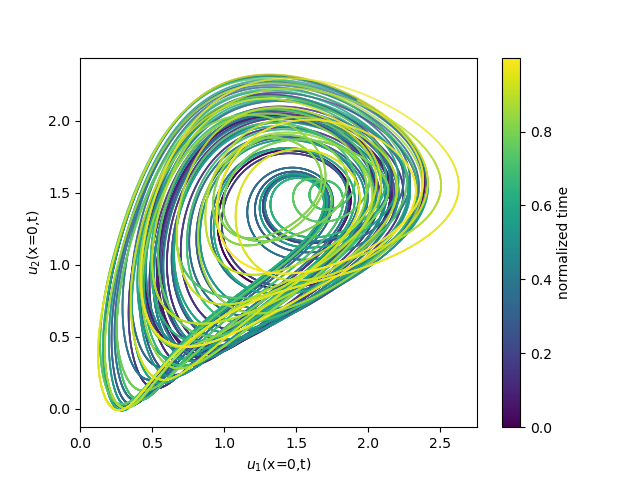}
\end{subfigure}
\begin{subfigure}{0.3\textwidth}
\includegraphics[width=\linewidth]{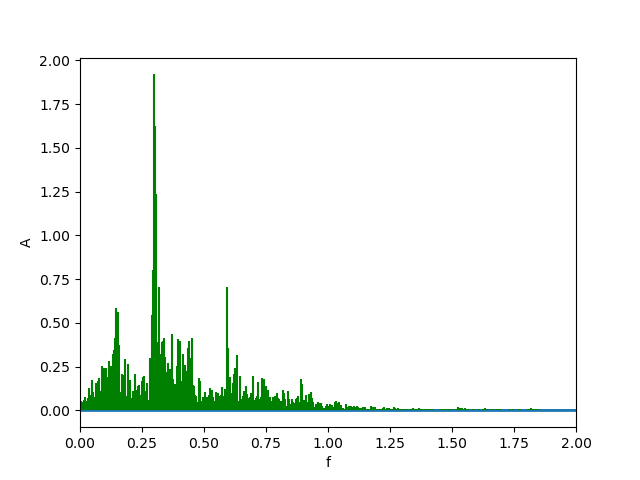}
\end{subfigure}
\begin{subfigure}{0.3\textwidth}
\includegraphics[width=\linewidth]{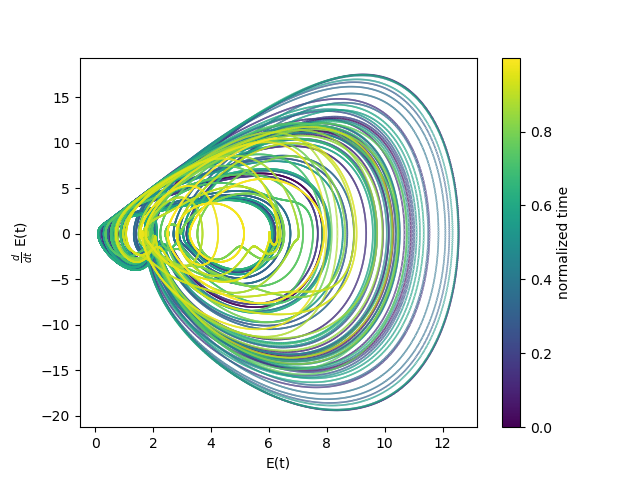}
\end{subfigure}
\caption{Chaos and strange attractors. 
Orbits of model (\ref{model}), starting from unstable periodic orbit of (Figure \ref{2cycle}),
in each of the four diffusion regimes, with $C = -1.5329, -1.53, -1.63, -1.51$.
For each orbit we display the local phase space (\ref{attractor_local}), the spectrum of energy (\ref{E}), 
and energy phase space (\ref{attractor_global}).
We observe the development of a strange attractor in every case.}
\label{fig:2Dattractors}
\end{figure}

Results are displayed on Figure \ref{fig:2Dattractors}.
In each of the four regimes of diffusion, we can observe similarities.
The mix of blue and yellow colors, on the local phase space, is an indication of orbit revisit.
Therefore orbits are bounded, and visited repeatedly through time.
The power spectrum is continuous (as opposed to clear separated lines), therefore the signal is not periodic.
Clear differences appear on the qualitative description of the attractor in energy space.

%%%%%%%%%%%%%%%%%%%%%%% 2D attractors %%%%%%%%%%%%%%%%%%%%%%%%%%%%%%%%%%%%%%%%

\section{Conclusion}

In this article, we have introduced an original model (\ref{model}) which contains a rich dynamics, 
exhibiting original Turing patterns, oscillating patterns, and chaos. 
We have derived accurate numerical methods in order to approximate orbits, steady states, 
and periodic orbits of the latter.
We also showed a convergence property (\ref{convergence}).
We have studied bifurcations of the model with respect to one control parameter, 
in four types of diffusion regimes: 
a linear case (that may be seen as a witness test), 
self-diffusion for each of the two species, and cross-diffusion between the species.
In each of the four cases, we have observed a three steps development.
Initially, a stable Turing pattern appears. 
Then, as the control parameter decreases, stability is lost through a Hopf bifurcation (Figure \ref{Hopf}). 
A periodic orbit appears, leading to an oscillating Turing pattern.
Eventually, as the control parameter decreases further, the periodic orbit loses its stability via a period doubling bifurcation (Figure \ref{1cycle}).
Decreasing further the control parameter, the double-period orbit itself loses stability (Figure \ref{2cycle}), but, interestingly, scenarios differ according to different regimes.
A period doubling cascade is observed (Figure \ref{road}) in the linear case, as well as in the first self-difusion case and the cross-diffusion case. To the best of author's knowledge, it is the first time that a clear raod to chaos is described, via bifurcation analysis, in the context of self- and cross-diffusion.
In the second self-diffusion case, a fold bifurcation is observed.
Eventually, a chaotic regime is reached in all cases, 
leading to the formation of chaotic attractors, differing according to the nonlinear regime that produced them.
During this study, we have used an original energetic approach, following \cite{Aymard:2022}. 
The latter has revealed to be more global than the local study, especially when studying the strange atttractors (Figure \ref{fig:2Dattractors}), allowing notably a much clearer separation in terms of strange attractor topology. 

As further works, several paths may be explored.
In the regime of self-diffusion of the second species, the road to chaos has not been clearly identified,
and therefore several options are possible, from the fold bifurcation point observed in Figure \ref{2cycle}.
Another interesting question, 
that arise naturally when comparing the strange attractors in Figure \ref{fig:2Dattractors},
would be to analyze their topology, using a branched manifold approach, 
in the spirit of \cite{Sciamarella:2020}, 
and find their classification in the taxonomy proposed in \cite{Letellier:2022}.

\end{document}